\documentstyle[amscd,amssymb,xypic,verbatim,11pt]{amsart}
\xyoption{all}

\newcommand{\nc}{\newcommand}
\newcommand{\rnc}{\renewcommand}

\nc{\exto}[1]{\stackrel{#1}{\longrightarrow}}
\nc{\dlim}{{\mathop{\lim\limits_{\longrightarrow}}}}
\nc{\lan}{\big\langle}
\nc{\ran}{\big\rangle}

\nc{\kk}{{\mathsf{k}}}
\nc{\ix}{{\mathsf{i}}}
\nc{\jx}{{\mathsf{j}}}

\nc{\C}{{\mathbb{C}}}
\nc{\HH}{{\mathbb{H}}}
\nc{\PP}{{\mathbb{P}}}
\nc{\QQ}{{\mathbb{Q}}}
\nc{\ZZ}{{\mathbb{Z}}}

\nc{\CA}{{\mathcal{A}}}
\nc{\CB}{{\mathcal{B}}}
\nc{\CC}{{\mathcal{C}}}
\nc{\D}{{\mathcal{D}}}
\nc{\CE}{{\mathcal{E}}}
\nc{\CF}{{\mathcal{F}}}
\nc{\CG}{{\mathcal{G}}}
\nc{\CH}{{\mathcal{H}}}
\nc{\CJ}{{\mathcal{J}}}
\nc{\CL}{{\mathcal{L}}}
\nc{\CM}{{\mathcal{M}}}
\nc{\CN}{{\mathcal{N}}}
\nc{\CO}{{\mathcal{O}}}
\nc{\CQ}{{\mathcal{Q}}}
\nc{\CR}{{\mathcal{R}}}
\nc{\CS}{{\mathcal{S}}}
\nc{\CT}{{\mathcal{T}}}
\nc{\CU}{{\mathcal{U}}}
\nc{\CV}{{\mathcal{V}}}
\nc{\CW}{{\mathcal{W}}}
\nc{\CX}{{\mathcal{X}}}
\nc{\CY}{{\mathcal{Y}}}
\nc{\CZ}{{\mathcal{Z}}}
\nc{\CMo}{{\mathcal{M}^\circ}}
\nc{\Co}{{{C}^\circ}}

\nc{\BY}{{\overline{Y}}}
\nc{\BZ}{{\overline{Z}}}
\nc{\BYD}{{\overline{Y}{}^{|D|}}}
\nc{\OZ}{{\overline{Z}}}
\nc{\bg}{{\bar{g}}}

\nc{\bq}{{\mathbf{q}}}
\nc{\BD}{{\mathbf{D}}}
\nc{\BG}{{\mathbf{G}}}
\nc{\BM}{{\mathbf{M}}}
\nc{\BP}{{\mathbf{P}}}
\nc{\BPr}{{\mathsf{P}}}
\nc{\BR}{{\mathbf{R}}}
\nc{\BRO}[1]{{{\mathbf{R}}^{\circ}_{#1}}}
\nc{\BRD}[1]{{{\mathbf{R}}^{|D|}_{#1}}}
\nc{\BRP}[1]{{{\mathbf{R}}^{1}_{#1}}}
\nc{\BRTP}[1]{{{\mathbf{\tilde{R}}}{}^{1}_{#1}}}
\nc{\BS}{{\mathbf{S}}}
\nc{\BMS}{{{\mathbf{M}}^{{s}}}}
\nc{\BMSS}{{{\mathbf{M}}^{{ss}}}}
\nc{\BMZ}{{\mathbf{M}^{\circ}}}
\nc{\BCL}{{\mathbf{L}}}

\nc{\PCC}{{{}^\perp\CC}}

\nc{\Cl}{{\mathsf{Cliff}}}
\nc{\Clev}{{\mathop{\mathsf{Cliff}}^{\circ}}}

\nc{\FA}{{\mathfrak{A}}}
\nc{\FB}{{\mathfrak{B}}}
\nc{\FI}{{\mathfrak{I}}}
\nc{\FZ}{{\mathfrak{Z}}}

\nc{\TFA}{{\tilde{\mathfrak{A}}}}
\nc{\TFB}{{\tilde{\mathfrak{B}}}}

\nc{\fa}{{\mathfrak{a}}}
\nc{\fg}{{\mathfrak{g}}}
\nc{\fp}{{\mathfrak{p}}}
\nc{\FD}{{\mathfrak{D}}}
\nc{\FE}{{\mathfrak{E}}}
\nc{\FL}{{\mathfrak{L}}}
\nc{\FM}{{\mathfrak{M}}}
\nc{\FR}{{\mathfrak{R}}}
\nc{\FS}{{\mathsf{S}}}

\nc{\sfc}{{\mathsf{c}}}
\nc{\sfch}{{\mathsf{ch}}}
\nc{\sfh}{{\mathsf{h}}}

\nc{\SK}{{\mathsf{K}}}
\nc{\SO}{{\mathsf{O}}}
\nc{\SQ}{{\mathsf{Q}}}
\nc{\SPV}{{\mathsf{S}^+\mathsf{V}}}
\nc{\SMV}{{\mathsf{S}^-\mathsf{V}}}
\nc{\SPMV}{{\mathsf{S}^\pm\mathsf{V}}}
\nc{\SX}{{S_X}}
\nc{\SY}{{S_Y}}
\nc{\phipsi}{{q}}
\nc{\eps}{\varepsilon}

\nc{\pim}{{\pi_-}}
\nc{\pip}{{\pi_+}}

\nc{\BE}{{\overline{\CE}}}
\nc{\TE}{{\tilde{\CE}}}
\nc{\TQ}{{\tilde{Q}}}
\nc{\TCF}{{\tilde{\CF}}}
\nc{\TCG}{{\tilde{\CG}}}
\nc{\TCL}{{\tilde{\CL}}}
\nc{\TF}{{\tilde{F}}}
\nc{\TW}{{\tilde{W}}}
\nc{\TCB}{{\widetilde{\CB}}}
\nc{\TCC}{{\tilde{\CC}}}
\nc{\TCX}{{\tilde{\CX}}}
\nc{\TCY}{{\tilde{\CY}}}
\nc{\TPhi}{{\tilde{\Phi}}}
\nc{\OPhi}{{\bar{\Phi}}}
\nc{\txi}{{\tilde{\xi}}}
\nc{\tp}{{\tilde{p}}}
\nc{\tq}{{\tilde{q}}}
\nc{\tzeta}{{\tilde{\zeta}}}
\nc{\tpi}{{\tilde{\pi}}}

\nc{\HCB}{{\widehat{\CB}}}
\nc{\HE}{{\widehat{\CE}}}
\nc{\HS}{{\widehat{S}}}
\nc{\HX}{{\hat{X}}}
\nc{\hxi}{{\hat{\xi}}}

\nc{\UH}{{\mathcal{H}}}

\nc{\TM}{{\widetilde{M}}}
\nc{\TCM}{{\widetilde{\CM}}}
\nc{\TS}{{\widetilde{S}}}
\nc{\TU}{{\widetilde{U}}}
\nc{\TX}{{\widetilde{X}}}
\nc{\TY}{{\widetilde{Y}}}
\nc{\TYO}{{{\widetilde{Y}}^\circ}}
\nc{\barf}{{\bar{f}}}
\nc{\te}{{\tilde{e}}{}}
\nc{\tf}{{\tilde{f}}}
\nc{\tg}{{\tilde{g}}}
\nc{\ti}{{\tilde{\imath}}}
\nc{\tj}{{\tilde{\jmath}}}
\nc{\ty}{{\tilde{y}}}
\nc{\tphi}{{\tilde{\phi}}}
\nc{\hf}{{\hat{f}}}

\nc{\urho}{{\underline{\rho}}}

\nc{\LRA}{\Leftrightarrow}
\nc{\RA}{\Rightarrow}
\nc{\lotimes}{\mathbin{\mathop{\otimes}\limits^{\mathbb{L}}}}
\nc{\CEnd}{\mathop{\mathcal{E}\mathit{nd}}\nolimits}
\nc{\CExt}{\mathop{\mathcal{E}\mathit{xt}}\nolimits}
\nc{\CHom}{\mathop{\mathcal{H}\mathit{om}}\nolimits}
\nc{\RH}{\mathop{{\mathsf{R}}\Gamma}\nolimits}
\nc{\RGamma}{\mathop{{\mathsf{R}}\Gamma}\nolimits}
\nc{\RHom}{\mathop{\mathsf{RHom}}\nolimits}
\nc{\RCHom}{\mathop{\mathsf{R}\mathcal{H}\mathit{om}}\nolimits}
\nc{\RG}{\mathop{\mathsf{R\Gamma}}\nolimits}
\nc{\Hom}{\mathop{\mathsf{Hom}}\nolimits}
\nc{\Ext}{\mathop{\mathsf{Ext}}\nolimits}
\nc{\End}{\mathop{\mathsf{End}}\nolimits}
\nc{\Tor}{\mathop{\mathsf{Tor}}\nolimits}
\nc{\Tordim}{\mathop{\mathsf{Tor}\text{\rm-}\mathsf{dim}}\nolimits}
\nc{\Hilb}{\mathop{\mathsf{Hilb}}\nolimits}
\nc{\Spec}{\mathop{\mathsf{Spec}}\nolimits}
\nc{\Proj}{\mathop{\mathsf{Proj}}\nolimits}
\nc{\Pic}{\mathop{\mathsf{Pic}}\nolimits}
\renewcommand{\Im}{\mathop{\mathsf{Im}}\nolimits}
\nc{\Tw}{\mathop{\mathsf{Tw}}\nolimits}
\nc{\Ker}{\mathop{\mathsf{Ker}}\nolimits}
\nc{\Coker}{\mathop{\mathsf{Coker}}\nolimits}
\nc{\codim}{\mathop{\mathsf{codim}}\nolimits}
\nc{\sing}{{\mathsf{sing}}}
\nc{\supp}{\mathop{\mathsf{supp}}}
\nc{\perf}{{\mathsf{perf}}}
\nc{\rank}{\mathop{\mathsf{rank}}}
\nc{\Pf}{{\mathsf{Pf}}}
\nc{\Gr}{{\mathsf{Gr}}}
\nc{\OGr}{{\mathsf{OGr}}}
\nc{\Flag}{{\mathsf{Fl}}}
\nc{\Kosz}{{\mathsf{Kosz}}}
\nc{\LGr}{{\mathsf{LGr}}}
\nc{\GTGr}{{\mathsf{G_2Gr}}}
\nc{\GTF}{{\mathsf{G_2F}}}
\nc{\OF}{{\mathsf{OF}}}
\nc{\Fl}{{\mathsf{Fl}}}
\nc{\Bl}{{\mathsf{Bl}}}
\nc{\GL}{{\mathsf{GL}}}
\nc{\PGL}{{\mathsf{PGL}}}
\nc{\SL}{{\mathsf{SL}}}
\nc{\SP}{{\mathsf{Sp}}}
\nc{\Spin}{{\mathsf{Spin}}}
\nc{\Tot}{{\mathsf{Tot}}}
\nc{\ev}{{\mathsf{ev}}}
\nc{\od}{{\mathsf{odd}}}
\nc{\coev}{{\mathsf{coev}}}
\nc{\id}{{\mathsf{id}}}
\nc{\opp}{{\mathsf{opp}}}
\nc{\PS}{{{\PP^3}}}
\nc{\Qu}{{{Q^3}}}
\nc{\tdim}{\mathop{\Tor\dim}}
\nc{\ecart}{{\fbox{$\scriptstyle\mathsf{EC}$}}}
\nc{\ad}{{\mathop{\mathsf ad}}}
\nc{\gr}{{\mathop{\mathsf gr}}}
\nc{\qgr}{{\mathop{\mathsf qgr}}}
\nc{\tor}{{\mathop{\mathsf tor}}}
\rnc{\mod}{{\mathop{\mathsf mod}}}
\nc{\Mod}{{\mathop{\mathsf Mod}}}
\nc{\Coh}{{\mathop{\mathsf Coh}}}
\nc{\Ab}{{\mathop{\mathcal{A}\mathit{b}}}}
\nc{\QCoh}{{\mathop{\mathsf QCoh}}}

\nc{\AAV}{{\mathcal{AAV}}}

\nc{\Rep}{{\mathsf{Rep}}}

\nc{\Cubics}{{{\mathcal{S}}_3}}
\nc{\VFT}{{{\mathcal{S}}_{14}}}
\nc{\VFTE}{{{\mathcal{N}}_{\mathrm{reg,sm}}}}
\nc{\MX}{{\CM_X}}
\nc{\MY}{{\CM_Y}}
\nc{\MYE}{{\CM_{Y,\CE}}}
\nc{\Yd}{{Y_d}}
\nc{\Yfive}{{Y_5}}
\nc{\Xg}{{X_{2g-2}}}
\nc{\Xtt}{{X_{22}}}
\nc{\Xst}{{X_{16}}}
\nc{\Xtw}{{X_{12}}}
\nc{\Xe}{{X_{8}}}
\nc{\Xf}{{X_{4}}}

\nc{\git}{{/\!\!/\!{}_\chi}}

\theoremstyle{plain}

\newtheorem{theorem}{Theorem}[section]

\newtheorem{lemma}[theorem]{Lemma}
\newtheorem{proposition}[theorem]{Proposition}
\newtheorem{corollary}[theorem]{Corollary}

\theoremstyle{definition}

\newtheorem{definition}[theorem]{Definition}

\newtheorem{example}[theorem]{Example}

\theoremstyle{remark}

\newtheorem{remark}[theorem]{Remark}

\newenvironment{proof}{\noindent{\sf Proof:}}{\qed\medskip}

\title[Derived categories of quadric fibrations and intersections of quadrics]%
{Derived categories of quadric fibrations\protect\\and intersections of quadrics}
\author{Alexander Kuznetsov}
\address{
Algebra Section, Steklov Mathematical Institute,
8 Gubkin str., Moscow 119991 Russia
}
\email{akuznet@@mi.ras.ru}
\date{}
\thanks{I was partially supported by RFFI grants 05-01-01034 and 02-01-01041,
Russian Presidential grant for young scientists No. MK-3926.2004.1,
CRDF Award No. RUM1-2661-MO-05, and the Russian Science Support Foundation.}

\begin{document}

\begin{abstract}
We construct a semiorthogonal decomposition of the derived category of coherent sheaves
on a quadric fibration consisting of several copies of the derived category of the base
of the fibration and the derived category of coherent sheaves of modules over the sheaf
of even parts of the Clifford algebras on the base corresponding to this quadric fibration
generalizing the Kapranov's description of the derived category of a single quadric.
As an application we verify that the noncommutative algebraic variety $(\PP(S^2W^*),\CB_0)$,
where $\CB_0$ is the universal sheaf of even parts of Clifford algebras,
is Homologically Projectively Dual to the projective space $\PP(W)$
in the double Veronese embedding $\PP(W) \to \PP(S^2W)$. Using the properties of
the Homological Projective Duality we obtain a description of the derived category
of coherent sheaves on a complete intersection of any number of quadrics.
\end{abstract}

\maketitle

\section{Introduction}

Derived categories of coherent sheaves on algebraic varieties recently
came into focus of renewed mathematical investigation. One of the reasons
for this is their role in the Homological Mirror Symmetry Conjecture
of M.~Kontsevich \cite{Ko} predicting an equivalence of the derived category
of coherent sheaves on an algebraic variety with the Fukaya category
of its mirror.

An important tool of investigation of the derived categories of coherent sheaves
is given by a notion of a semiorthogonal decomposition introduced in~\cite{B,BK}.
Semiorthogonal decomposition is a way to split a derived category into simpler
triangulated categories. On the other side of the mirror it corresponds to
a decomposition of the Fukaya category of a Landau--Ginzburg model
with respect to the singular fibers of a potential.

The first example of a nontrivial semiorthogonal decomposition
was discovered in~\cite{BO1}. It was shown there that the derived category
of a smooth complete intersection of two even-dimensional quadrics
contains the derived category of a hyperelliptic curve as a full subcategory
and the orthogonal is generated by an exceptional collection.
Further, in~\cite{BO2} Bondal and Orlov suggested a general description
of the derived category of coherent sheaves on a complete intersection
of any number of quadrics in terms of the sheaf of Clifford algebras
on the space of quadrics.

Other examples of semiorthogonal decompositions
for some Fano threefolds appeared in~\cite{K3,K4,K},
and a general method of constructing such decompositions
was given in~\cite{K2}. This method is based on a notion
of a {\sf Homological Projective Duality}\/ (HP-duality for short).
Algebraic varieties $X$ and $Y$ equipped with morphisms into dual
projective spaces $f:X \to \PP(V)$ and $g:Y \to \PP(V^*)$ are
{\sf HP-dual}\/ if the derived category of $Y$
can be embedded fully and faithfully into the derived category
of the universal hyperplane section $\CX \subset X\times\PP(V^*)$
of $X$ in a certain way. Properties of HP-dual varieties were thoroughly
investigated in~\cite{K2}. In particular, it was shown there that
for HP-dual varieties $X$ and $Y$, if
$X_L = X\times_{\PP(V)}\PP(L^\perp)$ and
$Y_L = Y\times_{\PP(V^*)}\PP(L)$ are orthogonal linear sections
of $X$ and $Y$ of expected dimensions
($L$ is a linear subspace in $V^*$ and $L^\perp \subset V$ is its orthogonal),
then there exist semiorthogonal decompositions of the derived categories
of $X_L$ and $Y_L$ with equivalent nontrivial terms.
The goal of this paper is to give a proof of the result of Bondal and Orlov
using the approach of HP-duality.

Explicitly, we take $X = \PP(W)$, a projective space,
$V = S^2W$ and $f:X \to \PP(V)$, the double Veronese embedding.
Then the universal hyperplane section of $\PP(W)$ with respect to the embedding $f$
is the universal quadric $\CX \subset \PP(W)\times\PP(S^2W^*)$, and $X_L$
is the intersection of quadrics parameterized by the space $L \subset S^2W^*$.
Further, we consider the sheaf of even parts of Clifford algebras on $\PP(S^2W^*)$
$$
\CB_0 =
\CO_{\PP(S^2W^*)} \oplus
\Lambda^2W \otimes \CO_{\PP(S^2W^*)}(-1) \oplus
\Lambda^4W \otimes \CO_{\PP(S^2W^*)}(-2) \oplus
\dots,
$$
prove that the noncommutative algebraic variety $Y = (\PP(S^2W^*),\CB_0)$
is HP-dual to $X$.

{\rnc{\thetheorem}{\ref{xyhpd}}
\begin{theorem}
The noncommutative algebraic variety $Y = (\PP(S^2W^*),\CB_0)$
is Homologically Projectively Dual to $\PP(W)$ in the double Veronese embedding.
\end{theorem}
}
Finally, we note that the semiorthogonal decompositions
of the intersection of quadrics $X_L$ and of $Y_L = (\PP(L),\CB_0)$
provided by the theory of HP-duality
coincide with those predicted by Bondal and Orlov.

{\rnc{\thetheorem}{\ref{dbxl}}
\begin{theorem}
For any vector subspace $L\subset S^2W^*$ such that the corresponding
intersection of quadrics $X_L$ is complete there exists
a semiorthogonal decomposition
$$
\begin{array}{llll}
\D^b(X_L) &=&  \langle \D^b(\PP(L),\CB_0),\CO_{X_L}(1),\dots,\CO_{X_L}(n-2r)\rangle, &
\text{if $r = \dim L \le n/2$}\\
\D^b(\PP(L),\CB_0) &=&  \langle \CB_{n-2r},\dots,\CB_{-2},\CB_{-1},\D^b(X_L)\rangle, &
\text{if $r = \dim L \ge n/2$}
\end{array}
$$
and an equivalence $\D^b(X_L) \cong \D^b(\PP(L),\CB_0)$ if $\dim L = n/2$,
where $\D^b(\PP(L),\CB_0)$ is the derived category of sheaves of $\CB_0$-modules on $\PP(L)$,
$\CB_1 =
W \otimes \CO_{\PP(S^2W^*)} \oplus
\Lambda^3W \otimes \CO_{\PP(S^2W^*)}(-1) \oplus
\Lambda^5W \otimes \CO_{\PP(S^2W^*)}(-2) \oplus
\dots$
is the sheaf of odd parts of Clifford algebras on $\PP(S^2W^*)$, and
$\CB_{k-2} = \CB_k\otimes\CO_{\PP(S^2W^*)}(-1)$ for all $k$.
\end{theorem}
}

Actually, in order to establish the HP-duality of $X = \PP(W)$ and $Y = (\PP(S^2W^*),\CB_0)$
we need to check that the universal quadric $\CX \subset \PP(W)\times\PP(S^2W^*)$
admits a certain semiorthogonal decomposition. We consider $\CX$ as a quadric fibration
over $\PP(S^2W^*)$ and instead of this particular case consider the general situation
of any algebraic variety $\CX$ equipped with a structure of a flat quadric fibration
$p:\CX \to S$ over any smooth base scheme $S$. With such a fibration we associate
a sheaf of algebras (the sheaf of even parts of Clifford algebras) $\CB_0$ on $S$
and show that the derived category of coherent sheaves of $\CB_0$-modules on $S$
gives a nontrivial semiorthogonal component in $\D^b(\CX)$.

{\rnc{\thetheorem}{\ref{dbcx}}
\begin{theorem}
If $p:\CX \to S$ is a flat quadric fibration of relative dimension $n-2$
then there exists a semiorthogonal decomposition
$$
\D^b(\CX) =
\lan
\D^b(S,\CB_0),
p^*(\D^b(S))\otimes\CO_{\CX/S}(1),
p^*(\D^b(S))\otimes\CO_{\CX/S}(2),
\dots,
p^*(\D^b(S))\otimes\CO_{\CX/S}(n-2)\ran,
$$
where $\D^b(S,\CB_0)$ is the derived category of coherent sheaves of $\CB_0$-modules on $S$.
\end{theorem}
}
The proof is based on the Koszul duality between the coordinate algebra
of $\CX$ over $S$ and the homogeneous Clifford algebra and follows
closely the Kapranov's description of the derived category
of a single quadric~\cite{Ka}.

The paper is organized as follows.
In section~2 we recall the necessary background.
In section~3 we discuss the Koszul duality between
the coordinate algebra of a quadric fibration
and the homogeneous Clifford algebra, and give
its geometric interpretation involving the sheaf
of even parts of the Clifford algebras.
In section~4 we give a description of the derived category
of coherent sheaves of a quadric fibration.
In section~5 we check that the noncommutative algebraic variety $(\PP(S^2W^*),\CB_0)$
is Homologically Projectively Dual to the projective space $\PP(W)$ in the double Veronese
embedding $\PP(W) \to \PP(S^2W)$ and deduce from this theorem of Bondal and Orlov.


{\bf Acknowledgements:}
I would like to thank A.~Bondal and D.~Orlov for useful discussions.


\section{Preliminaries}

Throughout the paper we assume the base field $\kk$ to be algebraically closed of characteristic $0$.

We start this section with a brief overview of what noncommutative algebraic geometry is.
The modern approach to the noncommutative algebraic geometry is the following:
we consider some abelian (or triangulated) category of noncommutative origin
and consider it as the category of coherent sheaves (or the derived category
of coherent sheaves) on the would-be noncommutative variety. Certainly
the category must be ``sufficiently nice''. There are several sources of
nice abelian categories. Let us briefly remind some of them.

\subsection{A sheaf of finite algebras}

Let $\BY$ be a usual algebraic variety, $\CB$ a sheaf of $\CO_\BY$-algebras
of finite rank as a $\CO_\BY$-module, and consider the category $\mod_\BY\text{-}\CB$
of coherent sheaves of right $\CB$-modules on $\BY$. This category is sufficiently nice
at least when $\CB$ is locally free over $\CO_\BY$. We consider it as the category
of coherent sheaves on a noncommutative algebraic variety $Y = (\BY,\CB)$:
$$
\Coh(\BY,\CB) = \mod_\BY\text{-}\CB.
$$
Besides the coherent category it is sometimes convenient to consider also the quasicoherent
category $\QCoh(\BY,\CB) = \Mod_\BY\text{-}\CB$ (we remove the condition of finite generatedness).
As usually (cf.~\cite{SGA6}) the derived category of coherent sheaves $\D^-(\Coh(\BY,\CB_\BY))$
can be identified with the full subacategory of $\D^-(\QCoh(\BY,\CB_\BY))$ consisting of objects
with coherent cohomologies. As we already mentioned above
noncommutative varieties of the form $(\BY,\CB)$ with $\CB$ being locally free over $\CO_\BY$
form a nice category. Let us briefly sketch some of the related definitions and constructions.

A {\sf morphism}\/ $f:(\BY,\CB_\BY) \to (\BZ,\CB_\BZ)$ is a pair $(f_\circ,f_\CB)$,
where $f_\circ:\BY\to \BZ$ is a morphism of algebraic varieties, and
$f_\CB:f_\circ^*\CB_\BZ \to \CB_\BY$ is a homomorphism of
$f_\circ^*\CO_\BZ\cong\CO_\BY$-algebras.
With every morphism $f:(\BY,\CB_\BY) \to (\BZ,\CB_\BZ)$ we can associate
the pushforward functor $R^0f_*:\QCoh(\BY,\CB_\BY) \to \QCoh(\BZ,\CB_\BZ)$ and
the pullback functor $L_0f^*:\QCoh(\BZ,\CB_\BZ) \to \QCoh(\BY,\CB_\BY)$,
where $\QCoh$ stands for the category quasicoherent sheaves as follows:
$$
R^0f_*(F) = R^0{f_\circ}_* F,
\qquad\text{and}\qquad
L_0f^*(G) = L_0f_\circ^*G \otimes_{f_\circ^*\CB_\BZ}\CB_\BY
$$
(the structures of a $\CB$-module on the RHS are evident).
It is clear that $R^0f_*$ is left exact and $L_0f^*$ is right exact.
On the other hand, the categories $\QCoh(\BY,\CB_\BY)$ and
$\QCoh(\BZ,\CB_\BZ)$ have enough injective and enough locally free objects
respectively, hence there exist the derived functors
$$
f_*:\D^+(\QCoh(\BY,\CB_\BY)) \to \D^+(\QCoh(\BZ,\CB_\BZ))
\qquad\text{and}\qquad
f^*:\D^-(\QCoh(\BZ,\CB_\BZ)) \to \D^-(\QCoh(\BY,\CB_\BY)).
$$
The pullback functor always takes
$\D^-(\Coh(\BZ,\CB_\BZ))$ to $\D^-(\Coh(\BY,\CB_\BY))$ and
when morphism $f$ is projective the pushforward functor takes
$\D^b(\Coh(\BY,\CB_\BY))$ to $\D^b(\Coh(\BZ,\CB_\BZ))$

In a contrast with commutative case one have to distinguish between the categories
$\Coh(\BY,\CB_\BY)$ and $\Coh(\BY,\CB_\BY^\opp)$ of right and left $\CB_\BY$-modules.
There are local $\RCHom$ and tensor product functors
$$
\arraycolsep=2pt
\begin{array}{rcl}
\RCHom_{\CB_\BY}  &:& \D^-(\Coh(\BY,\CB_\BY))\times\D^+(\Coh(\BY,\CB_\BY)) \to \D^+(\BY,\CO_\BY),\smallskip\\
\otimes_{\CB_\BY} &:& \D^-(\Coh(\BY,\CB_\BY))\times\D^-(\Coh(\BY,\CB_\BY^\opp)) \to \D^-(\BY,\CO_\BY),
\end{array}
$$
Also sometimes it is useful to consider the category $\Coh(\BY,\CB_\BY\otimes\CB_\BY^\opp)$
of $\CB_\BY$-bimodules. The functors $\RCHom_{\CB_\BY}$ and $\otimes_{\CB_\BY}$ can be defined
on the bimodule categories and take values in the appropriate module categories
(e.g. the tensor product over $\CB_\BY$ of two $\CB_\BY$-bimodules is again a $\CB_\BY$-bimodule).

Finally we should mention the twisted pullback functor
$f^!:\D^+(\Coh(\BZ,\CB_\BZ)) \to \D^+(\Coh(\BY,\CB_\BY))$.
It can be defined for embeddable morphisms along the lines of~\cite{H}.
When $f_\circ$ is smooth and $\CB_\BY = f_\circ^*\CB_\BZ$ we put
$f^!(G) = f^*G \otimes \omega_{\BY/\BZ}[\dim\BY/\BZ]$ and when
$f_\circ$ is finite we put
$$
f^!G = \RCHom_{\CB_\BZ}({f_\circ}_*\CB_\BY,G) \in
\D^+(\Coh(\BZ,{f_\circ}_*\CB_\BY)) = \D^+(\Coh(\BY,\CB_\BY)).
$$
Finally, arbitrary embeddable morphism can be decomposed into a composition
of a finite map and a smooth map and the twisted pullback can be defined
by functoriality. Standard verification shows that all usual relations between
functors still holds (see~\cite{H,K}).


\subsection{A sheaf of graded algebras}

Let $S$ be a usual algebraic variety, and $\FB = \oplus_{k=0}^\infty \FB_k$
a sheaf of graded $\CO_S$-algebras (with sheaves $\FB_k$ being locally free
over $\CO_S$ of finite rank). Consider the abelian category $\qgr_S\text{-}\FB$,
the quotient category of the category $\gr_S\text{-}\FB$ of finitely generated over $\CO_S$
sheaves of graded right $\FB$-modules on $S$ by the subcategory of sheaves of
graded $\FB$-modules of finite rank over $\CO_S$.
We consider it as the category of coherent sheaves on a noncommutative algebraic variety $\Proj_S(\FB)$:
$$
\Coh(\Proj_S(\FB)) = \qgr_S\text{-}\FB.
$$
This kind of categories was considered in~\cite{AZ} in case $S = \Spec\kk$.
It was shown there that they are sufficiently nice if the algebra $\FB$ is
strongly noetherian and satisfies so-called $\chi$-condition. A situation
is quite similar for sheaves of graded algebras over any $S$.
Let us briefly sketch some of the related definitions and constructions.

Let $\pi:\gr_S\text{-}\FB \to \qgr_S\text{-}\FB$ be the factorization functor.
First of all, there exists a projection from the noncommutative variety $\Proj_S(\FB)$
to the base $S$, $q:\Proj_S(\FB) \to S$. The corresponding pullback functor
$q^*:\Coh(S) \to \qgr_S\text{-}\FB$ is given by a simple formula
$$
q^*(F) = \pi(F\otimes_{\CO_S}\FB).
$$
Note that in our assumptions on $\FB$ the functor $q^*$ is exact,
so we can consider the projection $q$ to be flat.

The pushforward functor $q_*:\qgr_S\text{-}\FB \to \Coh(S)$ is given by the formula
$$
\Gamma(U,q_*(\pi(M))) = \Hom_{\qgr_U\text{-}\FB_U}(\FB_U,M_U),
$$
where $U$ is an affine open subset of $S$, $M$ is a sheaf of graded $\FB$-modules on $S$,
and $\FB_U$, $M_U$ denote the restrictions of $\FB$, $M$ to $U$.

Similarly, we define a $S$-local $\Hom$-functor
$\CHom_{S}:(\qgr_S\text{-}\FB)^\opp\times\qgr_S\text{-}\FB \to \Coh(S)$ as follows
$$
\Gamma(U,\CHom_{S}(\pi(M),\pi(N))) = \Hom_{\qgr_U\text{-}\FB_U}(M_U,N_U),
$$
and a local $\Hom$-functor
$\CHom_Y:(\qgr_S\text{-}\FB)^\opp\times\qgr_S\text{-}\FB \to \qgr_S\text{-}\FB$ as follows
$$
\CHom_Y(M,N) = \pi\left(\oplus_{k=0}^\infty \CHom_{S}(M,N\lan k\ran)\right),
$$
where $\lan k\ran$ is the autoequivalence of $\qgr_S\text{-}\FB$
induced by the shift of grading functor on $\gr_S\text{-}\FB$.
Then
$$
\CHom_{S}(-,-) = q_*\CHom_Y(-,-),
\qquad\text{and}\qquad
\Hom_{\qgr_S\text{-}\FB}(-,-) = \Gamma(S,\CHom_{S}(-,-)).
$$

\begin{definition}[cf.~\cite{AZ}]
A sheaf of $\CO_S$-algebras $\FB$ is called {\sf strongly noetherian}\/
if for any base change $T \to S$ the sheaf of $\CO_T$-algebras $\FB_T$
(the pullback of $\FB$) is noetherian.

A sheaf of $\CO_S$-algebras $\FB$ is said to satisfy the {\sf relative $\chi$-condition}\/
if $\Ext^i_{\mod\text{-}\FB_U}(\CO_U,M_U)$ is a finitely generated $\CO_U$-module
for all affine $U\subset S$, all $i\ge 0$, and all $M\in\mod\text{-}\FB$.
\end{definition}

It is straightforward to check that all arguments of~\cite{AZ} work for the category $\qgr_S\text{-}\FB$
if the sheaf of algebras $\FB$ is strongly noetherian and satisfies the relative $\chi$-condition.

It is worth mentioning that if a sheaf of algebras is commutative and generated by its
first component (linearly generated) then the category $\qgr_S\text{-}\FB$ is equivalent
to the category of coherent sheaves on the usual (commutative) $\Proj_S(\FB)$ scheme
(this is the famous Theorem of Serre). However, even for commutative but not linearly
generated sheaf of algebras we can get an interesting example of a noncommutative variety.

\begin{example}\label{stack}
Let $\CL$ and $\CM$ be line bundles on $S$ and $d$ a global section of the bundle
$\CM^{*2}\otimes\CL^{*n}$. Let
$$
\FB = S^\bullet(\CL^*)\otimes(\CO_S \oplus \CM)
$$
with the structure of an algebra given by the map $d:\CM^2 \to \CL^{*n}$
and with the grading induced by $\deg\CL^* = 2$, $\deg\CM = n$.
It is easy to check that locally (over $S$) the category $\qgr_S\text{-}\FB$
can be identified with the category of coherent sheaves on the double covering
of $S$ ramified in the zero locus $Z$ of $d$ $\ZZ/2\ZZ$-equivariant
with respect to the involution of the double covering (note that when $n$ is odd
this covering in general doesn't exist globally).
Thus $\qgr_S\text{-}\FB$ is equivalent locally to the category of coherent sheaves
on the quotient stack of the double covering by the involution.
We will call this noncommutative variety ``$S$ with a $\ZZ/2\ZZ$-stack structure along $Z$''.
\end{example}


%

Finally we have to say some words about the double version of the category $\qgr$.
Let $\FB$ be a sheaf of bigraded $\CO_S$-algebras. Consider the category $\gr_S^2\text{-}\FB$
of right bigraded sheaves of $\FB$-modules and its subcategory $\tor_S^2\text{-}\FB$
consisting of all sheaves of bigraded $\FB$-modules $M$ such that $M_{k,l} = 0$
for $k\ge k_0$, $l\ge l_0$ where $k_0$ and $l_0$ are sufficiently large.
Let (cf.\ \cite{BGK})
$$
\qgr_S^2\text{-}\FB := \gr_S^2\text{-}\FB / \tor_S^2\text{-}\FB
$$
be the quotient category. It is shown in~\cite{BGK} that under appropriate
conditions category $\qgr^2$ behaves analogously to $\qgr$.
Moreover, it is shown there that there exists an equivalence of categories
$$
\qgr_S^2\text{-}\FB \cong \qgr_S\text{-}(\Delta(\FB)),
$$
where $\Delta(\FB) = \oplus_{k=0}^\infty \FB_{k,k}$ is the diagonal subalgebra of $\FB$.

\subsection{Koszul Duality for sheaves of graded algebras}

Let $S$ be a scheme and $\FA = \FA_0 \oplus \FA_1 \oplus \FA_2 \oplus \dots$ be a sheaf
of graded $\CO_S$-algebras.
$\FA$ is called {\sf flat of finite type}\/ if all components $\FA_k$ of $\FA$
are locally free of finite rank as $\CO_S$-modules.
$\FA$ is called {\sf connected}\/ if $\FA_0 = \CO_S$.
$\FA$ is called {\sf linearly generated}\/ if it is connected and the canonical
homomorphism of algebras $T^\bullet_{\CO_S}(\FA_1) \to \FA$ from the tensor algebra
of $\FA_1$ over $\CO_S$ to $\FA$ is surjective.
$\FA$ is called {\sf quadratic}\/ if it is connected, linearly generated,
and the kernel of the canonical homomorphism $T^\bullet_{\CO_S}(\FA_1) \to \FA$
is generated by its degree $2$ component
$$
\FI_\FA^2 =
T^2_{\CO_S}(\FA_1) \cap \Ker(T^\bullet_{\CO_S}(\FA_1) \to \FA)  =
\Ker(\FA_1\otimes_{\CO_S} \FA_1 \to \FA_2) \subset
\FA_1\otimes_{\CO_S} \FA_1.
$$
Let $\FA$ be a quadratic flat sheaf of $\CO_S$-algebras of finite type.
The {\sf quadratic dual}\/ sheaf of $\CO_S$-algebras $\FA^!$ of $\FA$
is defined as the quotient of the tensor algebra of $\FA_1^*$ by the ideal
generated by the subbundle
$$
(\FI_A^2)^\perp :=
\Ker((\FA_1\otimes_{\CO_S} \FA_1)^* \to (\FI_\FA^2)^*) =
\FA_2^* \subset
(\FA_1\otimes_{\CO_S} \FA_1)^* \cong
\FA_1^*\otimes_{\CO_S} \FA_1^*.
$$
Thus
$$
\FA^! := T^\bullet_{\CO_S}(\FA_1^*)/\lan (\FI_\FA^2)^\perp \ran.
$$
Let $\FA$ be a quadratic flat sheaf of $\CO_S$-algebras of finite type
and let $\FA^!$ be its quadratic dual.
The multiplications in algebras $\FA^!$ and $\FA$ give maps
$\FA^!_k\otimes\FA_1^* \to \FA^!_{k+1}$ and $\FA_k\otimes\FA^{!*}_1 \to \FA_{k+1}$.
Dualizing, we obtain maps
$\FA^{!*}_{k+1} \to (\FA^!_k)^*\otimes\FA_1$ and $\FA_{k+1}^* \to \FA_k^*\otimes\FA^!_1$,
which induce morphisms of graded modules
$\FA^{!*}_{k+1}\otimes\FA \to \FA^{!*}_k\otimes\FA\lan1\ran$ and
$\FA_{k+1}^*\otimes\FA^! \to \FA_k^*\otimes\FA^!\lan1\ran$,
where $\langle k\rangle$ denotes the shift of grading by~$k$.
On the other hand, we have the augmentation maps
$\FA^{!*}_0\otimes\FA \cong \FA \to \FA_0 = \CO_S$ and
$\FA_0^*\otimes\FA^! \cong \FA^! \to \FA^!_0 = \CO_S$
which are morphisms of graded modules as well.
Obtained in this way sequences of $\FA$-modules and of $\FA^!$-modules
\begin{equation}\label{koszul}
\begin{array}{r}
\dots \to
\FA^{!*}_3\otimes\FA\langle-3\rangle \to
\FA^{!*}_2\otimes\FA\langle-2\rangle \to
\FA^{!*}_1\otimes\FA\langle-1\rangle \to
\FA^{!*}_0\otimes\FA \to
\CO_S \to 0,\\
\dots \to
\FA_3^*\otimes\FA^!\langle-3\rangle \to
\FA_2^*\otimes\FA^!\langle-2\rangle \to
\FA_1^*\otimes\FA^!\langle-1\rangle \to
\FA_0^*\otimes\FA^! \to
\CO_S \to 0,
\end{array}
\end{equation}
are complexes, called the {\sf Koszul complexes}\/ of $\FA$ and $\FA^!$ respectively.

\begin{definition}[cf.~\cite{PP}]
A quadratic flat sheaf of $\CO_S$-algebras of finite type $\FA$ is called
{\sf Koszul}\/ if the Koszul complexes~$(\ref{koszul})$ are acyclic.
\end{definition}

Denote by $\FR^k_\FA$ the cohomologies of the
truncated in degree $k$ and shifted Koszul complexes of $\FA$
$$
\FR^k_\FA =
\Ker(\FA^{!*}_{k}\otimes\FA \to \FA^{!*}_{k-1}\otimes\FA\langle1\rangle) =
\Coker(\FA^{!*}_{k+2}\otimes\FA\langle-2\rangle \to \FA^{!*}_{k+1}\otimes\FA\langle-1\rangle).
$$
Then we have the following left and right resolutions
\begin{equation}\label{leftr}
\dots \to
\FA^{!*}_{k+3}\otimes\FA\langle-3\rangle \to
\FA^{!*}_{k+2}\otimes\FA\langle-2\rangle \to
\FA^{!*}_{k+1}\otimes\FA\langle-1\rangle \to
\FR^k_\FA \to 0,
\end{equation}
\begin{equation}\label{rightr}
0 \to \FR^k_\FA \to
\FA^{!*}_{k  }\otimes\FA \to
\FA^{!*}_{k-1}\otimes\FA\langle1\rangle \to
\dots \to
\FA^{!*}_{0  }\otimes\FA\langle k\rangle \to
\CO_S\lan k\ran \to 0.
\end{equation}

The importance of modules $\FR^k_\FA$ is demonstrated
by the following relative generalization of the result of~\cite{KKO}.
Consider the category of bigraded $\FA$-bimodules on $S$.
Consider the {\sf diagonal}\/ object
$$
\FA_\Delta = \oplus_{p,q\ge 0} \FA^{p+q}
$$
in this category. As we will see below the diagonal object plays the role
of the structure sheaf of the diagonal on the square of a variety and
the resolution of the following theorem gives a resolution of the diagonal.

\begin{theorem}\label{resdiag}
The diagonal object $\FA_\Delta$ admits the following
resolution in the category of bigraded $\FA$-bimodules on $S$:
\begin{equation}\label{rfa}
\dots \to
\FR^2_\FA \otimes \FA\langle-2\rangle \to
\FR^1_\FA \otimes \FA\langle-1\rangle \to
\FR^0_\FA \otimes \FA \to
\FA_\Delta \to 0,
\end{equation}
\end{theorem}
\begin{proof}
This is proved in \cite{KKO}, propositions 4.7 and 4.9
for Koszul algebras over a field. The same arguments work
as well for Koszul flat sheaves of algebras of finite type.
\end{proof}

\subsection{Clifford algebras}\label{clalg}

Let $V$ be a $\kk$-vector space and $q\in S^2V^*$, a quadratic form on $V$.
The Clifford algebra of $q$ is defined as follows (see~\cite{Bou})
$$
B_q = T^\bullet(V)/\langle v\otimes v' + v'\otimes v - 2q(v,v')\cdot 1 \rangle_{v,v'\in V},
$$
where $T^\bullet(V)$ is a free associative algebras generated by $V$
and $1$ is the unit of $T^\bullet(V)$.
The Clifford algebra $B_q$ is a finite dimensional associative $\kk$-algebra
(if $v_1,\dots,v_n$ is a basis of $V$ then the monomials
$v_{i_1}v_{i_2}\dots v_{i_r}$, $1\le i_1 < i_2 < \dots < i_r \le n$,
form a basis of $B_q$ over $\kk$).
The Clifford algebra is naturally $\ZZ/2\ZZ$-graded:
its even component $B_{q0}$ is spanned by monomials $v_{i_1}v_{i_2}\dots v_{i_r}$ of even degree,
and its odd component $B_{q1}$ is spanned by monomials of odd degree.

If the quadratic form $q$ is nondegenerate then both the Clifford algebra $B_q$
and its even part $B_{q0}$ are semisimple~\cite{Bou}.
If $n = 2m$ then $B_q$ is isomorphic to a matrix algebra of rank $2^m$ and
$B_{q0}$ is isomorphic to the product of two matrix algebras of rank $2^{m-1}$.
Similarly, if $n = 2m+1$ then $B_q$ is isomorphic to a product of two matrix algebras
of rank $2^m$ and $B_{q0}$ is isomorphic to a matrix algebra of rank $2^m$.

On a contrary, if $q = 0$ is a zero quadratic form, then it is clear that the Clifford algebra
is isomorphic to the exterior algebra $\Lambda^\bullet(V)$, so has infinite homological dimension.

\begin{lemma}\label{multv}
Let $v\in V$ and consider {\rm(}either left or right{\rm)} multiplication
by $v$ maps $B_{q0} \stackrel{v}\to B_{q1}$ and $B_{q1} \stackrel{v}\to B_{q0}$.
If $q(v) \ne 0$ then both are isomorphisms and
if $q(v) = 0$ then both have rank $2^{n-2}$.
\end{lemma}
\begin{proof}
The square of these maps is equal to the multiplication by $q(v)$, hence for $q(v)\ne 0$
the maps are isomorphisms. Now assume that $q(v) = 0$ and consider for instance the left
multiplication by $v$. Choose a basis $v_1,\dots,v_n$ of $V$ such that $v_1 = v$.
Then it is clear that both the kernel and the image of the multiplication maps
are spanned by monomials $v_1v_{i_2}\dots v_{i_r}$ with $2\le i_2 < \dots < i_r\le n$,
hence are $2^{n-2}$-dimensional.
\end{proof}

Now let $q$ be any quadratic form. Let $v_1,\dots,v_n$ be an orthogonal for $q$ basis of $V$,
so that we have
$$
v_iv_j = - v_jv_i,
\quad\text{for $i\ne j$}
\quad\text{and}\quad
v_i^2 = q(v_i)
\qquad\text{in $B_q$.}
$$
Let
$$
d := v_1v_2\cdots v_n \in B_q.
$$
Then it is easy to see that $v_i d = (-1)^{n-1} d v_i$.
Therefore $d$ is central in $B_q$ for odd $n$
and central in $B_{q0}$ for even $n$.
Moreover, it is easy to check that $d$ (up to a scalar) doesn't depend
on a choice of the orthogonal basis.
Note also that $d^2 = (-1)^{n(n-1)/2}\det q$.

\begin{lemma}\label{cldeg}
Assume that $q$ has $1$-dimensional kernel. Then
for $n=2m$ the quotient algebra $B_{q0}/B_{q0}d$ is isomorphic to a matrix algebra of rank $2^{m-1}$ and
for $n=2m+1$ the quotient algebra $B_{q}/B_{q}d$ is isomorphic to a matrix algebra of rank $2^{m}$.
\end{lemma}
\begin{proof}
Choose an orthogonal basis $v_1,\dots,v_n$ in $V$ such that $v_n$ generates the kernel of $q$.
It is clear that the ideal $B_{q0}d$ (resp.\ $B_{q}d$) is spanned by
monomials $v_{i_1}v_{i_2}\dots v_{i_{r-1}}v_n$ with $1\le i_1 < i_2 < \dots < i_{r-1} < n$.
Therefore the quotient algebra $B_{q0}/B_{q0}d$ (resp.\ $B_{q}/B_{q}d$) is spanned by
monomials $v_{i_1}v_{i_2}\dots v_{i_{r}}$ with $1\le i_1 < i_2 < \dots < i_{r} \le n$,
hence isomorphic to (the even part of) the Clifford algebra of the restriction of $q$
to the vector subspace of $V$ spanned by $v_1,\dots,v_{n-1}$.
But the restriction of $q$ to this vector subspace is nondegenerate,
hence the corresponding algebra is indeed a matrix algebra.
\end{proof}

\section{Clifford algebra of a quadric fibration}

We need to introduce some notation.
Let
\begin{itemize}
\item $S$ be a smooth algebraic variety;
\item $E$ be a vector bundle of rank $n$ on $S$;
\item $\pi:\PP_S(E) \to S$ be the projectivization of $E$ on $S$;
\item $\CO_{\PP_S(E)/S}(1)$ be the Grothendieck line bundle on $\PP_S(E)$;
\item $\CL$ be a line bundle on $S$;
\item $\sigma:\CL \to S^2E^*$ be an embedding of vector bundles\\(in the other words
$\sigma \in \Gamma(S,S^2E^*\otimes \CL^*) = \Gamma(\PP_S(E),\CO_{\PP_S(E)/S}(2)\otimes \CL^*)$);
\item $\CX \subset \PP_S(E)$ be the zero locus of $\sigma$ on $\PP_S(E)$; and
\item $p:\CX \to S$ be the restriction of $\pi$ to $\CX$.
\end{itemize}
The projection $p:\CX \to S$ is a {\sf quadric fibration}.
Note that $p$ is a flat morphism since $\sigma:\CL \to S^2E^*$ has no zeroes.
The relative dimension of $p$ equals $n-2$.


\subsection{Koszul duality for the coordinate algebra of $\CX$ over $S$}\label{kd}

Consider the sheaf of graded algebras
$$
\FA_\sigma =
\bigoplus_{k=0}^\infty \FA_{\sigma k} =
\bigoplus_{k=0}^\infty p_*(\CO_{\CX/S}(k)) \cong
\bigoplus_{k=0}^\infty S^kE^*/(S^{k-2}E^*\otimes \sigma(\CL)) \cong
T^\bullet E^* /
\langle \Lambda^2E^* \oplus \sigma(\CL) \rangle,
$$
where $T^\bullet E^* = \bigoplus_{k=0}^\infty {E^*}^{\otimes k}$
is the sheaf of free associative algebras generated over $\CO_S$ by $E^*$.
It is clear that $\FA_\sigma$ is a quadratic flat sheaf of $\CO_S$-algebras of finite type.
Consider the sheaf of its quadratic dual algebras
$$
\FB_\sigma = \FA_\sigma^! =
T^\bullet E /
\langle \Ker(\sigma^*:S^2E \to \CL^*) \rangle.
$$
We call $\FB_\sigma$ {\sf the sheaf of homogeneous Clifford algebras of $\sigma$}\/
(cf.~\cite{Ka}). Note that as a sheaf of $\CO_S$-modules $\FB_{\sigma k}$ takes form
\begin{equation}\label{fb_k}
\FB_{\sigma k} \cong \Lambda^kE \oplus \Lambda^{k-2}E \otimes \CL^* \oplus \Lambda^{k-4}E\otimes \CL^{*2} \oplus \dots,
\end{equation}
so that
\begin{equation}\label{fb_s}
\FB_\sigma \cong \CO_S \oplus E \oplus (\Lambda^2E \oplus \CL^*) \oplus
(\Lambda^3E \oplus E \otimes \CL^*) \oplus
(\Lambda^4E \oplus \Lambda^2E \otimes \CL^* \oplus \CL^{*2}) \oplus \dots
\end{equation}
with a Clifford multiplication.

\begin{lemma}
The sheaves of $\CO_S$-algebras $\FA_\sigma$ and $\FB_\sigma$ on $S$ are Koszul.
\end{lemma}
\begin{proof}
See e.g.~\cite{PP}~II.6.
\end{proof}

Further we will sometimes omit the index $\sigma$ to unburden the notation.

\begin{lemma}
The algebras $\FA$ and $\FB$ are strongly noetherian and satisfy
the relative $\chi$-condition.
\end{lemma}
\begin{proof}
Under any change of base $T \to S$ the algebras $\FA$ and $\FB$
go to the algebras of the same type corresponding to the quadric fibration $\CX_T \to T$, so it suffices to check the usual
noetherian property for $\FA$ and $\FB$. For this we note
that the algebra $\FA$ is commutative and finitely generated
(over $\CO_S$), hence noetherian. On the other hand, the algebra
$\FB$ is finitely generated as a module over its central subalgebra
$S^\bullet(\CL^*) \subset \FB$. The algebra $S^\bullet(\CL^*)$ is noetherian
by the same reasons as $\FA$ is, hence $\FB$ is noetherian as well.
To check the relative $\chi$-condition we just use the Koszul complexes
for computing $\Ext$'s in the category of $\FB = \FA^!$-modules,
considering them as free resolutions of $\CO_S$.
\end{proof}

Consider the global section $\delta$ of the vector bundle
$E^*\otimes E = \FA_1\otimes\FB_1$ on $S$ corresponding to the identity endomorphism of $E$.
Consider a structure of an algebra on $S^\bullet(E^*)\otimes\FB$ in which
$S^\bullet(E^*)$ and $\FB$ commute.

\begin{lemma}\label{delta2}
We have $\delta\cdot\delta = \sigma \in S^2E^*\otimes \CL^* \subset S^2E^*\otimes \FB_2$.
\end{lemma}
\begin{proof}
Let $e_i$ be a basis of $E$ and let $x^i$ be the dual basis of $E^*$. Then we have
$$
\delta\cdot\delta =
(\sum_i x^i\otimes e_i)(\sum_j x^j\otimes e_j) =
\sum_{i,j}x^ix^j\otimes e_ie_j =
\frac12\sum_{i,j}x^ix^j\otimes (e_ie_j + e_je_i) =
\sum_{i,j}\sigma(e_i,e_j)x^ix^j =
\sigma.
$$
since $x^ix^j = x^jx^i$ and $e_ie_j + e_je_i = \sigma(e_i,e_j)$.
\end{proof}

Let $\delta$ and $\delta'$ denote the multiplication by $\delta$ from the left
and from the right in the algebra $S^\bullet(E^*)\otimes\FB$ as homomorphisms of right (resp.\ left)
bigraded $S^\bullet(E^*)\otimes\FB$-modules on $S$:
$$
\delta,\delta': S^\bullet(E^*)\otimes\FB \to S^\bullet(E^*)\otimes\FB\lan1,1\ran.
$$
It follows from the above lemma that
\begin{equation}\label{dlr2}
\delta\circ\delta = \sigma,\qquad
\delta'\circ\delta' = \sigma
\end{equation}
Note also that $\sigma$ is {\em central}\/ in $S^\bullet(E^*)\otimes\FB$.

Now we are going to give a geometric meaning to the categories
of graded $\FA$ and $\FB$-modules and to the corresponding Koszul complexes.
For a start consider the category of $\FA$-modules.


\subsection{Geometric interpretation for $\FA$}


The Theorem of Serre for the sheaf of algebras $\FA$ takes the following form:

\begin{proposition}
There is an exact equivalence of abelian categories
$$
\qgr_S\text{-}\FA \cong \Coh(\CX)
$$
taking a sheaf $F$ on $\CX$ to the graded $\FA$-module
$\bigoplus_{k=0}^\infty p_*(F\otimes\CO_{\CX/S}(k))$.
In particular, the line bundles $\CO_{\CX/S}(k)$ correspond
to the shifted free module $\FA\lan k\ran$, and moreover,
the twist functor $F\mapsto F\otimes\CO_{\CX/S}(k)$ in the category $\Coh(\CX)$
corresponds to the shift of grading functor $M \mapsto M\lan k\ran$
in the category $\qgr_S\text{-}\FA$.
\end{proposition}


Similarly, we can consider the category of bigraded $\FA$-bimodules.


\begin{proposition}
There is an exact equivalence of abelian categories
$$
\qgr^2_S\text{-}(\FA^\opp\otimes_{\CO_S}\FA) \cong \Coh(\CX\times_S\CX)
$$
taking a sheaf $F$ on $\CX\times_S\CX$ to the bigraded $\FA$-bimodule
$\bigoplus_{k,l=0}^\infty (p\times p)_*(F\otimes\CO_{\CX\times_S\CX/S}(k,l))$.
In particular, the exterior product $F\boxtimes G$ of sheaves $F,G\in\Coh(\CX)$
corresponds to the tensor product {\rm(}over $\CO_S${\rm)} of the corresponding $\FA$-modules,
and the structure sheaf of the diagonal $\Delta_*\CO_\CX$
corresponds to the diagonal object~$\FA_\Delta$.
\end{proposition}

Let $\CR^k_\FA$ denote the sheaf on $\CX$ corresponding
to the graded $\FA$-module $\FR^k_\FA$ defined by~(\ref{leftr}) or (\ref{rightr}).
Applying the equivalences of the above propositions
to the exact sequences~(\ref{leftr}), (\ref{rightr}) and~(\ref{rfa})
we obtain the following exact sequences on $\CX$
\begin{equation}\label{cral}
\dots \to
\FB_{k+3}^*\otimes\CO_{\CX/S}(-3) \to
\FB_{k+2}^*\otimes\CO_{\CX/S}(-2) \to
\FB_{k+1}^*\otimes\CO_{\CX/S}(-1) \to
\CR^k_\FA \to 0,
\end{equation}
\begin{equation}\label{crar}
0 \to \CR^k_\FA \to
\FB_{k}^*\otimes\CO_\CX \to
\FB_{k-1}^*\otimes\CO_{\CX/S}(1) \to
\dots \to
\FB_{0}^*\otimes\CO_{\CX/S}(k) \to 0,
\end{equation}
and on $\CX\times_S\CX$
\begin{equation}\label{resda}
\dots \to
\CR^2_\FA \boxtimes \CO_{\CX/S}(-2) \to
\CR^1_\FA \boxtimes \CO_{\CX/S}(-1) \to
\CR^0_\FA \boxtimes \CO_\CX \to
\Delta_*\CO_\CX \to 0.
\end{equation}

\subsection{Geometric interpretation for $\FB$}\label{gfb}


We start with the following important observation helps

\begin{lemma}
The subalgebra $S^\bullet(\CL^*) = \CO_S \oplus \CL^* \oplus \CL^{*2} \oplus \dots \subset \FB$
is central.
\end{lemma}
\begin{proof}
A direct calculation in the homogeneous Clifford algebra.
\end{proof}

Now we can consider $\FB$ as a $S^\bullet(\CL^*)$-module.
It turns out that $\FB$ is finitely generated over $S^\bullet(\CL^*)$.
Geometrically this means that noncommutative variety $\Proj_S(\FB)$ is
a finite covering of $\Proj_S(S^\bullet(\CL^*)) = S$. Thus, the category
$\qgr_S\text{-}\FB$ can be identified with the category of coherent sheaves
of modules over a certain sheaf of algebras on $S$.

Explicitly, consider the pushforward functor $q_*:\qgr_S\text{-}\FB \to \Coh(S)$.
It is isomorphic to the composition of the restriction functor
$\qgr_S\text{-}\FB \to \qgr_S\text{-}S^\bullet(\CL^*)$
with the equivalence $\qgr\text{-}S^\bullet(\CL^*) \to \Coh(S)$,
and can be written as
\begin{equation}\label{gpf}
M = \oplus M_k \mapsto \dlim (M_{2k}\otimes \CL^k),
\end{equation}
where the inductive limit is considered with respect to the maps
$M_{2k}\otimes \CL^k \to M_{2(k+1)}\otimes \CL^{k+1}$
induced by the $S^\bullet(\CL^*)$-module structure on $M$ (note that
$\CL^* \subset \FB$ sits in the degree 2).
Under this functor the free module $\FB$ goes to $\dlim (\FB_{2k}\otimes \CL^k)$
which according to (\ref{fb_k}) is isomorphic to
\begin{equation}\label{cb0}
\CB_0 := \CO_S \oplus (\Lambda^2E \otimes \CL) \oplus (\Lambda^4E \otimes \CL^2) \oplus \dots,
\end{equation}
{\sf the even part of the Clifford algebra of $\sigma$}.
It is clear that the pushforward of any graded $\FB$-module has a structure
of a $\CB_0$-module on $S$, so we can consider the functor $q_*$ as a functor
$q_*:\qgr_S\text{-}\FB \to \mod_S\text{-}\CB_0$.

\begin{proposition}\label{gs}
The functor $q_*:\qgr_S\text{-}\FB \to \mod_S\text{-}\CB_0$, $M \mapsto \dlim (M_{2k}\otimes \CL^k)$
is an exact equivalence of abelian categories. Furthermore, we have
\begin{equation}\label{qsml}
q_*(M\lan 2l\ran) \cong q_*(M) \otimes \CL^{*l}
\quad\text{for all $M \in \qgr_S\text{-}\FB$, $l\in \ZZ$,}
\end{equation}
$$
q_*\CHom_{\qgr_S\text{-}\FB}(M,N) \cong \CHom_{\CB_0}(q_*M,q_*N),
\quad\text{for all $M,N \in \qgr_S\text{-}\FB$.}
$$
\end{proposition}
\begin{proof}
It is convenient to decompose the functor $q_*$ as a composition of the functor
$$
\qgr_S\text{-}\FB \to \qgr_S\text{-}\FB^{(2)},\qquad \oplus M_k \mapsto \oplus M_{2k},
$$
where $\FB^{(2)} = \oplus_{k=0}^\infty \FB_{2k}$ is the double Veronese
subalgebra of $\FB$, with the functor
$$
\alpha:\qgr_S\text{-}\FB^{(2)} \to \mod\text{-}\CB_0,\qquad
\oplus M_k \mapsto \dlim (M_k\otimes \CL^k).
$$
The first functor is an exact equivalence by \cite{P}, so it remains
to consider the second. To this end we consider also a functor
$$
\beta:\mod\text{-}\CB_0 \to \qgr_S\text{-}\FB^{(2)},\qquad
\CM \mapsto \oplus_{k=0}^\infty (\CM\otimes \CL^{*k}).
$$
The structure of a graded module over the Veronese subalgebra
$\FB^{(2)} = \oplus_{k=0}^\infty \FB_{2k}$ on $\oplus_{k=0}^\infty (\CM \otimes \CL^{*k})$
is defined by the maps
$$
(\CM\otimes \CL^{*k}) \otimes \FB_{2l} \cong
(\CM\otimes \CL^{*(k+l)}) \otimes (\FB_{2l} \otimes \CL^l) \to
(\CM\otimes \CL^{*(k+l)}) \otimes \CB_0 \to
\CM\otimes \CL^{*(k+l)},
$$
where the second map is the canonical morphism
$\FB_{2l} \otimes \CL^l \to \dlim(\FB_{2k} \otimes \CL^k) = \CB_0$
(in our case it is just an embedding
$\CO_S \oplus (\Lambda^2E \otimes \CL) \oplus \dots \oplus (\Lambda^{2k}E \otimes \CL^k) \subset
\CO_S \oplus (\Lambda^2E \otimes \CL) \oplus \dots \oplus (\Lambda^{2k}E \otimes \CL^k) \oplus \dots$),
and the third morphism is given by the action of $\CB_0$ on $\CM$.

It is clear that $\alpha(\beta(\CM)) = \dlim (\CM\otimes \CL^{*k}\otimes \CL^k) = \dlim\, \CM = \CM$, so
$$
\alpha\circ\beta \cong \id.
$$
It remains to check that $\beta\circ\alpha \cong \id$.
To this end we note that the canonical maps $M_k\otimes \CL^k \to \dlim (M_l\otimes \CL^l)$
induce maps $M_k \to (\dlim (M_l\otimes \CL^l)) \otimes \CL^{*k}$ for all $k$
which together give a morphism of graded $\FB^{(2)}$-modules $\phi_M:M \to \beta(\alpha(M))$.
Applying the functor $\alpha$ to the exact sequence
$$
0 \to \Ker \phi_M \to M \to \beta(\alpha(M)) \to \Coker\phi_M \to 0
$$
we deduce that $\alpha(\Ker \phi_M) = \alpha(\Coker \phi_M) = 0$ since $\alpha$ is exact
and $\alpha\circ\beta \cong \id$ as we have noticed above. Thus it remains to check that
$\alpha(M) = 0$ implies $M = 0$. Indeed, if $\alpha(M) = \dlim (M_k\otimes \CL^k) = 0$
then the finite generatedness of $M$ implies that $M_k = 0$ for $k \gg 0$, hence $M = 0$
in $\qgr_S\text{-}\FB^{(2)}$.

Further, we have
$q_*(M\lan 2l\ran) :=
\dlim (M_{2k+2l}\otimes \CL^k) \cong
(\dlim (M_{2k+2l}\otimes \CL^{k+l})) \otimes \CL^{*l} \cong
q_*(M) \otimes \CL^{*l}$.

Finally, for affine $U\subset S$ we have
$\Gamma(U,q_*\CHom_{\qgr_S\text{-}\FB}(M,N)) =
\dlim(\Hom_{\qgr_U\text{-}\FB_U}(M_U,N_U\langle 2k\rangle)\otimes\CL_U^k)$,
but
\begin{multline*}
\Hom_{\qgr_U\text{-}\FB_U}(M_U,N_U\langle 2k\rangle) \cong
\Hom_{\CB_{0U}}(q_*(M_U),q_*(N_U\langle 2k\rangle)) \cong
\\ \cong
\Hom_{\CB_{0U}}(q_*(M_U),q_*(N_U)\otimes \CL_U^{*k}) \cong
\Hom_{\CB_{0U}}(q_*(M_U),q_*(N_U))\otimes \CL_U^{*k},
\end{multline*}
whereof we get
$\Gamma(U,q_*\CHom_{\qgr_S\text{-}\FB}(M,N)) =
\Hom_{\CB_{0U}}(q_*(M_U),q_*(N_U))$
for any affine $U$, and we obtain the desired isomorphism.
\end{proof}

Consider the graded $\FB$-module $\FB\lan s\ran$
($\lan s\ran$ is the shift of grading by $s$).
Under the functor $q_*$ it goes to the $\CB_0$-module
$$
\CB_s := q*(\FB\lan s\ran) = \dlim (\FB_{2k+s}\otimes \CL^k).
$$
For $s=1$ it is easy to deduce from (\ref{fb_k}) that
\begin{equation}\label{cb1}
\CB_1 \cong E \oplus (\Lambda^3E \otimes \CL) \oplus (\Lambda^5E \otimes \CL^2) \oplus \dots,
\end{equation}
{\sf the odd part of the Clifford algebra of $\sigma$}.
It follows from~(\ref{qsml}) that
\begin{equation}\label{gsok}
\CB_k \cong \begin{cases}
\CB_0\otimes \CL^{*l}, & \text{if $k = 2l$}\\
\CB_1\otimes \CL^{*l}, & \text{if $k = 2l + 1$}
\end{cases}
\end{equation}

\begin{lemma}\label{rhomb}
All sheaves\/ $\CB_k$ on $S$ are locally projective over $\CB_0$ and we have
$$
\RCHom_{\CB_0}(\CB_k,\CB_l) \cong \CB_{l-k}.
\qquad\text{for all $k,l\in\ZZ$}.
$$
\end{lemma}
\begin{proof}
We have $\CHom_{\CB_0}(\CB_k,q_*(-)) \cong q_*\CHom_{\qgr_S\text{-}\FB}(\FB\lan k\ran,-)$
by proposition~\ref{gs}. On the other hand, it is clear that
$\CHom_{\qgr_S\text{-}\FB}(\FB\lan k\ran,F) \cong F\langle-k\rangle$,
so it is an exact functor. Thus $\CHom_{\CB_0}(\CB_k,-)$ is an exact functor,
hence $\CB_k$ is locally projective over $\CB_0$.
Using proposition~\ref{gs} again we deduce
$$
\RCHom_{\CB_0}(\CB_k,\CB_l) \cong
q_*\CHom_{\qgr_S\text{-}\FB}(\FB\lan k\ran,\FB\lan l\ran) \cong
q_*(\FB\lan l-k\ran) \cong
\CB_{l-k}.
$$
for all $k,l\in\ZZ$.
\end{proof}

\begin{corollary}\label{tensorb}
All sheaves\/ $\CB_k$ on $S$ are flat over $\CB_0$ and we have
$$
\CB_k \otimes_{\CB_0} \CB_l \cong \CB_{k+l}
\qquad\text{for all $k,l\in\ZZ$}.
$$
\end{corollary}
\begin{proof}
We have
$\CB_k\otimes_{\CB_0} F \cong
\RCHom_{\CB_0}(\CB_{-k},\CB_0) \otimes_{\CB_0} F \cong
\RCHom_{\CB_0}(\CB_{-k},F)$
since $\CB_{-k}$ is locally projective.
The RHS functor is exact hence $\CB_k$ is flat over $\CB_0$.
Moreover, for $F = \CB_l$ the RHS equals $\CB_{k+l}$ by the above lemma.
\end{proof}

Similarly, we can consider the category of bigraded $\FB$-bimodules.


\begin{proposition}
There is an exact equivalence of abelian categories
$$
\qgr_S^2\text{-}(\FB^\opp\otimes\FB) \cong \mod\text{-}(\CB_0^\opp\otimes\CB_0)
$$
taking a bigraded bimodule $M = \oplus M_{i,j}$ to the sheaf
$(q\times q)_*(M) = \dlim (M_{2k,2k}\otimes \CL^{2k})$ which is naturally
a $\CB_0$-bimodule on $S$.
In particular, the exterior tensor product $F\boxtimes G$ of
$\CB_0$-modules $F\in\mod\text{-}\CB_0^\opp$,
$G\in\mod\text{-}\CB_0$ corresponds to the tensor product {\rm(}over $\CO_S${\rm)}
of the corresponding $\FB$-modules, and the diagonal bimodule $\CB_0$ corresponds
to the diagonal object $\FB_\Delta$.
\end{proposition}
\begin{proof}
Similarly.
\end{proof}

Let $\CR^k_\FB$ denote the $\CB_0$-module on $S$ corresponding
to the graded $\FB$-module $\FR^k_\FB$ defined by~(\ref{leftr}) or (\ref{rightr}).
Applying the equivalences of the above propositions
to the exact sequences~(\ref{leftr}), (\ref{rightr}) and~(\ref{rfa})
we obtain the following exact sequences of $\CB_0$-modules
\begin{equation}\label{crbl}
\dots \to
\FA_{k+3}^*\otimes\CB_{-3} \to
\FA_{k+2}^*\otimes\CB_{-2} \to
\FA_{k+1}^*\otimes\CB_{-1} \to
\CR^k_\FB \to 0,
\end{equation}
\begin{equation}\label{crbr}
0 \to \CR^k_\FB \to
\FA_{k}^*\otimes\CB_0 \to
\FA_{k-1}^*\otimes\CB_{1} \to
\dots \to
\FA_{0}^*\otimes\CB_{k} \to 0,
\end{equation}
and of $\CB_0$-bimodules
\begin{equation}\label{resdb}
\dots \to
\CR^2_\FB \boxtimes \CB_{-2} \to
\CR^1_\FB \boxtimes \CB_{-1} \to
\CR^0_\FB \boxtimes \CB_0 \to
\CB_0 \to 0.
\end{equation}

\subsection{Geometric interpretation for $\FA\otimes\FB$}\label{fafb}

Finally, we consider the categories of bigraded $S^\bullet(E^*)\otimes\FB$-modules
and $\FA\otimes\FB$-modules on $S$.

\begin{proposition}\label{ab_mod}
There is an exact equivalence of abelian categories
$$
\begin{array}{llllrrr}
\qgr_S^2\text{-}(S^\bullet(E^*)\otimes\FB) &\cong& \Coh(\PP_S(E),\CB_0),
\quad&\text{such that}&\quad
S^\bullet(E^*)\otimes\FB\lan k,l\ran &\mapsto& \CO_{\PP_S(E)/S}(k)\otimes\CB_l,\\
\qgr_S^2\text{-}(\FA\otimes\FB) &\cong& \Coh(\CX,\CB_0),
\quad&\text{such that}&\quad
\FA\otimes\FB\lan k,l\ran &\mapsto& \CO_{\CX/S}(k)\otimes\CB_l.
\end{array}
$$
\end{proposition}
\begin{proof}
Similarly to proposition~\ref{gs}.
\end{proof}

Recall the morphisms $\delta,\delta':S^\bullet(E)\otimes\FB \to S^\bullet(E)\otimes\FB\langle1,1\rangle$
constructed in section~\ref{kd}. They give rise to morphisms
$\delta_{k,l},\delta'_{k,l}:\CO_{\PP_S(E)/S}(k-1)\otimes\CB_{l-1} \to \CO_{\PP_S(E)/S}(k)\otimes\CB_{l}$
of right and left $\CB_0$-modules on $\PP_S(E)$.
Then lemma~\ref{delta2} implies the following

\begin{corollary}\label{d2s}
We have morphisms
$\delta_{k,l},\delta'_{k,l}:\CO_{\PP_S(E)/S}(k-1)\otimes\CB_{l-1} \to \CO_{\PP_S(E)/S}(k)\otimes\CB_{l}$
of right and left $\CB_0$-modules on $\PP_S(E)$ such that the compositions
$$
\delta_{k+1,l+1}\circ\delta_{k,l},
\delta'_{k+1,l+1}\circ\delta'_{k,l} :
\CO_{\PP_S(E)/S}(k-1)\otimes\CB_{l-1} \to \CO_{\PP_S(E)/S}(k+1)\otimes\CB_{l+1}
$$
coincide with the map given by
$\sigma \in \Gamma(\PP_S(E),\CO_{\PP_S(E)/S}(2)\otimes \CL^*)$.
\end{corollary}

\subsection{Additional central reduction: even-dimensional case}\label{ld}

Let us describe more explicitly the noncommutative algebraic variety $(S,\CB_0)$
in some special cases. Let $S_d \subset S$ denote the $d$-th degeneration locus of $\sigma$.
This is a closed subscheme of $S$ defined by the sheaf of ideals
$$
\CJ_d = \Im(
\xymatrix@1{\Lambda^{n+1-d}E\otimes\Lambda^{n+1-d}E\otimes \CL^{n+1-d} \ar[rr]^-{\Lambda^{n+1-d}\sigma} && \CO_S}
).
$$
In particular, $S_1$ is the zero locus of $\det\sigma\in\Gamma(S,\CL^{*n}\otimes(\det E^*)^2)$.

Assume that $n = \rank E = 2m$ is even.
Then it is easy to see that the sheaf of algebras
$$
\CB_0 =
\CO_S \oplus
\Lambda^2E \otimes \CL \oplus
\Lambda^4E \otimes \CL^2 \oplus
\dots \oplus
\Lambda^{2m}E \otimes \CL^m
$$
has a central subalgebra
$\CZ = \CO_S \oplus \det E \otimes \CL^m$
(locally it is generated by $1$ and the element $d$ defined in section~\ref{clalg}).
Therefore, we can identify $\CB_0$ with a pushforward
of a certain sheaf of algebras $\TCB_0$ on the double covering
$\phi:\TS = \Spec_{S}(\CZ) \to S$ ramified over $S_1$ (note that $d^2 = \det\sigma$).
In particular, the category $\Coh(S,\CB_0)$
identifies with the category $\Coh(\TS,\TCB_0)$.

Let $\TS_2\subset \TS$ be the preimage of the locus of quadrics of corank $\ge 2$
under the double covering $\phi:\TS \to S$, $\TS_2 = \phi^{-1}(S_2)$.

\begin{proposition}
The restriction of the sheaf of algebras $\TCB_0$ to $\TS\setminus\TS_2$
is a sheaf of Azumaya algebras.
\end{proposition}
\begin{proof}
By~\cite{M} it suffices to check that the fiber of the sheaf $\TCB_0$ at any point
$x \in \TS\setminus\TS_2$ is a matrix algebra of rank $2^{m-1}$.
Let $s = \phi(x) \in S \setminus S_2$. First, assume that $s\not\in S_1$, i.e.\
that $\sigma(s) \in S^2E_s^*$ is a nondegenerate quadric.
Then as it was mentioned in section~\ref{clalg} the even part of the corresponding Clifford algebra
$\CB_{\sigma(s) 0}$ is a product of two matrix algebras of rank $2^{m-1}$.
It is clear that the fiber of $\TCB_0$ at $x$ is just one of these matrix algebras.
On the other hand, assume that $s\in S_1\setminus S_2$, i.e.\
that $\sigma(s) \in S^2E_s^*$ is a quadric of corank $1$.
Since the fiber of the sheaf $\CZ$ over $s$ is spanned by $1$ and $d$,
the fiber of the algebra $\TCB_0$ at $x$ is the quotient
$\CB_{\sigma(s) 0}/\CB_{\sigma(s) 0}d$ which by lemma~\ref{cldeg} is also isomorphic
to a matrix algebra of rank $2^{m-1}$.
\end{proof}

Let us say that a quadric fibration $\CX \to S$ {\sf has only simple degenerations}\/
if all fibers are quadrics of corank $\le 1$ and the zero locus of $\det\sigma$ is smooth
(so that $\TS$ is also smooth).

If $p:\CX \to S$ has only simple degenerations then $\TCB_0$ is a sheaf
of Azumaya algebras on a smooth algebraic variety $\TS$. If $\TCB_0$ splits
then we have $\Coh(\TS,\TCB_0) \cong \Coh(\TS)$. This is true for example
when $\dim S = 1$ (\cite{M}).

\begin{corollary}\label{dc}
Assume that $S$ is smooth, $\dim S = 1$ and the quadric fibration
$\CX \to S$ has only simple degenerations. Then
$\Coh(S,\CB_0) \cong \Coh(\TS,\TCB_0) \cong \Coh(\TS)$.
\end{corollary}

\subsection{Additional central reduction: odd-dimensional case}


Assume that $n = \rank E = 2m+1$ is odd.
Then it is easy to see that the center
of the homogeneous Clifford algebra $\FB$ contains a subalgebra
$$
\FZ = S^\bullet(\CL^*) \otimes (\CO_S \oplus \det E).
$$
The corresponding category $\qgr_S\text{-}\FZ$ can be considered
as a category of coherent sheaves on $S$ with a $\ZZ/2\ZZ$-stack structure
along $S_1$ (see example~\ref{stack}). Let $\HS$ denote the corresponding
noncommutative variety (so that $\Coh(\HS) = \qgr_S\text{-}\FZ$)
and let $\psi:\HS \to S$ be the projection.
The sheaf of algebras $\FB$ considered as a $\FZ$-module gives rise to a sheaf
of algebras $\HCB_0$ on $\HS$ such that $\psi_*(\HCB_0) = \CB_0$. In particular,
the category $\Coh(S,\CB_0)$ identifies with the category $\Coh(\HS,\HCB_0)$.
Let $\HS_2\subset \HS$ be the preimage of the locus of quadrics of corank $\ge 2$
under the map $\psi:\HS \to S$, $\HS_2 = \psi^{-1}(S_2)$.

\begin{proposition}
The restriction of the sheaf of algebras $\HCB_0$ to $\HS\setminus\HS_2$
is a sheaf of Azumaya algebras.
\end{proposition}
\begin{proof}
The claim is local in $S$ so we may assume that $S$ is affine.
Then there exists a double covering $\phi:\TS \to S$ ramified in $S_1$
and $\HS$ is just the quotient stack of $\TS$ with respect to the action
of the group $\ZZ/2\ZZ$ generated by the involution of $\TS$ over $S$.
In particular, the projection $\phi$ factors as $\xymatrix@1{\TS \ar[r]^{\eta} & \HS \ar[r]^\psi & S}$,
where $\eta$ is \'etale. Consider the sheaf of algebras $\eta^*\HCB_0$ on $\TS$.
It suffices to check that $\eta^*\HCB_0$ is a sheaf of Azumaya algebras.

Note that the category $\Coh(\HS)$ is the category of $\ZZ/2\ZZ$-equivariant sheaves
on $\TS$, the functor $\psi_*$ boils down to the taking of invariants
with respect to the $\ZZ/2\ZZ$ action, and the functor $\eta^*$
is the forgetting of the $\ZZ/2\ZZ$-action.

Trivializing the sheaf $\CL$ on $S$ (and shrinking $S$ if necessarily) we can consider
the whole sheaf of Clifford algebras $\CB_\sigma$ of $\sigma$ on $S$.
It has a central subalgebra $\CZ = \CO_S \oplus \det E$
(locally it is generated by $1$ and the element $d$ defined in section~\ref{clalg})
and it is clear that $\Spec_{S}(\CZ) = \TS$. Therefore, we can identify
$\CB_\sigma$ with a pushforward of a certain sheaf of algebras $\TCB_\sigma$ on $\TS$.
This sheaf is naturally $\ZZ/2\ZZ$-graded, hence carries a $\ZZ/2\ZZ$-action and
can be considered as a sheaf of algebras on $\HS$. Its pushforward to $S$
(the invariants of $\CB_\sigma$ with respect to the $\ZZ/2\ZZ$-action) coincides
with the even part $\CB_{\sigma0}$ of the sheaf of Clifford algebras. This shows
that $\HCB_0$ coincides with $\TCB_\sigma$ considered as a $\ZZ/2\ZZ$-equivariant
sheaf of algebras on $\TS$. Hence $\eta^*\HCB_0 = \TCB_\sigma$ and it remains to check
that $\TCB_\sigma$ is a sheaf of Azumaya algebras on $\TS$.
This is done completely analogous with the even-dimensional case.

By~\cite{M} it suffices to check that the fiber of the sheaf $\CB_\sigma$ at any point
$x \in \TS\setminus\TS_2$ is a matrix algebra of rank $2^{m}$.
Let $s = \phi(x) \in S \setminus S_2$. First, assume that $s\not\in S_1$, i.e.\
that $\sigma(s) \in S^2E_s^*$ is a nondegenerate quadric. Then as it was mentioned
in section~\ref{clalg} the corresponding Clifford algebra $\CB_{\sigma(s)}$ is
a product of two matrix algebras of rank $2^{m}$. It is clear that the fiber of
$\TCB_\sigma$ at $x$ is just one of these matrix algebras.
On the other hand, assume that $s\in S_1\setminus S_2$, i.e.\
that $\sigma(s) \in S^2E_s^*$ is a quadric of corank $1$.
Since the fiber of the sheaf $\CZ$ over $s$ is spanned by $1$ and $d$,
the fiber of the algebra $\TCB_\sigma$ at $x$ is the quotient
$\CB_{\sigma(s)}/\CB_{\sigma(s)}d$ which by lemma~\ref{cldeg} is also isomorphic
to a matrix algebra of rank $2^{m}$.
\end{proof}

Similarly to corollary~\ref{dc} we deduce

\begin{corollary}\label{sc}
Assume that $S$ is smooth, $\dim S = 1$ and the quadric fibration
$\CX \to S$ has only simple degenerations. Then $\Coh(S,\CB_0) \cong \Coh(\HS)$.
\end{corollary}

\section{Derived category of a quadric fibration}


Recall the notation. Let
$S$ is a smooth algebraic variety;
$E$ is a vector bundle of rank $n$ on $S$;
$\pi:\PP_S(E) \to S$ is the projectivization of $E$ on $S$;
$\CO_{\PP_S(E)/S}(1)$ is the Grothendieck line bundle on $\PP_S(E)$;
$\CL$ is a line bundle on $S$;
$\sigma:\CL \to S^2E^*$ is an embedding of vector bundles;
$\CX \subset \PP_S(E)$ is the zero locus of $\sigma$ on $\PP_S(E)$; and
$p:\CX \to S$ is the restriction of $\pi$ to $\CX$.
We also denote by
\begin{itemize}
\item $i:\CX \to \PP_S(E)$ the embedding.
\end{itemize}
Then we have $p = \pi\circ i$.


The goal of this section is to describe a semiorthogonal decomposition
of the derived category of coherent sheaves on $\CX$, $\D^b(\CX)$.
\begin{definition}[\cite{B,BK}]
A semiorthogonal decomposition of a triangulated category $\CT$
is a sequence of full subcategories
$\CA_1,\dots,\CA_n$ in $\CT$ such that $\Hom_{\CT}(\CA_i,\CA_j) = 0$ for $i > j$
and for every object $T \in \CT$ there exists a chain of morphisms
$0 = T_n \to T_{n-1} \to \dots \to T_1 \to T_0 = T$ such that
the cone of the morphism $T_k \to T_{k-1}$ is contained in $\CA_k$
for each $k=1,2,\dots,n$.
\end{definition}
If a sequence $\CA_1,\dots,\CA_n$ forms a semiorthogonal decomposition of $\CT$
we write $\CT = \lan\CA_1,\CA_2,\dots,\CA_n\ran$.

Recall a sheaf of Clifford algebras $\CB_0$ on $S$ constructed in the previous section.

\begin{theorem}\label{dbcx}
If $p:\CX \to S$ is a flat quadric fibration of relative dimension $n-2$
then there exists a semiorthogonal decomposition
$$
\D^b(\CX) =
\lan
\D^b(S,\CB_0),
p^*(\D^b(S))\otimes\CO_{\CX/S}(1),
p^*(\D^b(S))\otimes\CO_{\CX/S}(2),
\dots,
p^*(\D^b(S))\otimes\CO_{\CX/S}(n-2)\ran,
$$
where $\D^b(S,\CB_0)$ is the derived category of coherent sheaves of $\CB_0$-modules on $S$.
\end{theorem}

The proof of the theorem takes the rest of the section.
We start with the following.

\begin{lemma}\label{poxs}
We have
$$
p_*(\CO_{\CX/S}(m)) =
\begin{cases}
\FA_m, & \text{for $m\ge 0$}\\
0, & \text{for $3-n \le m \le -1$}\\
\FA_{2-m-n}^* \otimes \det E \otimes \CL [2-n], & \text{for $m \le 2-n$}
\end{cases}
$$
\end{lemma}
\begin{proof}
We apply the functor $\pi_*$ to the exact sequence
$$
0 \to \CO_{\PP_S(E)/S}(m-2)\otimes \CL \to \CO_{\PP_S(E)/S}(m) \to i_*\CO_{\CX/S}(m) \to 0
$$
and take into account an isomorphism $p_*\CO_{\CX/S}(m) \cong \pi_*i_*\CO_{\CX/S}(m)$.
For $m\ge 0$ we get a distinguished triangle
$$
S^{m-2}E^*\otimes \CL \to S^mE^* \to p_*(\CO_{\CX/S}(m))
$$
which gives the first line in the formula;
for $3-n \le m \le -1$ we got the second line of the formula;
and for $m\le 2-n$ we get a distinguished triangle
$$
S^{2-m-n}E \otimes \det E \otimes \CL[1-n] \to S^{-m-n}E \otimes \det E[1-n] \to p_*(\CO_{\CX/S}(m)),
$$
which gives the third line of the formula.
\end{proof}

\begin{lemma}\label{easy}
The functor $\D^b(S) \to \D^b(\CX)$, $F \mapsto p^*F \otimes\CO_{\CX/S}(k)$,
is fully faithful for any $k\in \ZZ$ and the collection of subcategories
\begin{equation}\label{rm2col}
\hspace{-10pt}
\lan p^*(\D^b(S))\otimes\CO_{\CX/S}(k+1),p^*(\D^b(S))\otimes\CO_{\CX/S}(k+2),\dots,
p^*(\D^b(S))\otimes\CO_{\CX/S}(k+n-2)\ran
\end{equation}
is semiorthogonal in $\D^b(\CX)$.
\end{lemma}
\begin{proof}
Note that
$$
\Hom_\CX(p^*F \otimes\CO_{\CX/S}(k),p^*G \otimes\CO_{\CX/S}(l)) \cong
\Hom_\CX(p^*F,p^*G \otimes\CO_{\CX/S}(l-k)) \cong
\Hom_S(F,p_*(p^*G \otimes\CO_{\CX/S}(l-k))).
$$
By projection formula we have
$p_*(p^*G \otimes\CO_{\CX/S}(l-k)) \cong G \otimes p_*(\CO_{\CX/S}(l-k))$.
Since by lemma~\ref{poxs} we have $p_*(\CO_{\CX/S}(l-k)) = 0$ for $k-n+2 < l < k$,
the semiorthogonality follows;
and since $p_*(\CO_{\CX}) = \CO_S$ we deduce that
$\Hom_\CX(p^*F \otimes\CO_{\CX/S}(k),p^*G \otimes\CO_{\CX/S}(k)) \cong \Hom_S(F,G)$,
i.e. that the full faithfulness of the functor $F \mapsto p^*F\otimes\CO_{\CX/S}(k)$.
\end{proof}

Now we are going to construct a fully faithful functor
$\D^b(S,\CB_0) \to \D^b(\CX)$.
Recall the morphisms of right and left $\CB_0$-modules on $\PP_S(E)$
constructed in section~\ref{fafb}:
$$
\delta_{k,l},\delta'_{k,l}:\CO_{\PP_S(E)/S}(k-1)\otimes \CB_{l-1} \to \CO_{\PP_S(E)/S}(k)\otimes\CB_{l}.
$$

\begin{lemma}\label{kcd}
We have $\Ker \delta_{k,l} = \Ker \delta'_{k,l} = 0$, and
the sheaves $\Coker \delta_{k,l}$ and $\Coker \delta'_{k,l}$ are supported scheme-theoretically
on the hypersurface $\CX \subset \PP_S(E)$,
$$
\Coker \delta_{k,l} \cong i_*\CE_{k,l},
\qquad
\Coker \delta'_{k,l} \cong i_*\CE'_{k,l},
$$
where $\CE_{k,l}$ is a sheaf of right $\CB_0$-modules on $\CX$,
and $\CE'_{k,l}$ is a sheaf of left $\CB_0$-modules on $\CX$.
Moreover on $\CX$ there exist exact sequences of $\CB_0$-modules
$$
0 \to \CE_{k-1,l-1} \to \CO_{\CX/S}(k)\otimes\CB_{l} \to \CE_{k,l} \to 0,
\qquad
0 \to \CE'_{k-1,l-1} \to \CO_{\CX/S}(k)\otimes\CB_{l} \to \CE'_{k,l} \to 0
$$
and the periodicity isomorphisms
\begin{equation}\label{period}
\CE_{k,l}\otimes(\CO_{\CX/S}(s)\otimes \CL^t) \cong \CE_{k+s,l-2t},
\qquad
\CE'_{k,l}\otimes(\CO_{\CX/S}(s)\otimes \CL^t) \cong \CE'_{k+s,l-2t}.
\end{equation}
\end{lemma}
\begin{proof}
By corollary~\ref{d2s} we have the following commutative diagram
$$
\xymatrix{
0 \to \Ker \delta_{k,l} \ar[r] \ar[d]^\sigma &
\CO_{\PP_S(E)/S}(k-1)\otimes \CB_{l-1} \ar[r]^{\delta_{k,l}} \ar[d]^\sigma &
\CO_{\PP_S(E)/S}(k)\otimes \CB_{l} \ar[r] \ar[d]^\sigma \ar[dl]_{\delta_{k+1,l+1}} &
\Coker \delta_{k,l} \to 0 \ar[d]^\sigma \\
0 \to \Ker \delta_{k+2,l+2} \ar[r]  &
\CO_{\PP_S(E)/S}(k+1)\otimes \CB_{l+1} \ar[r]^{\delta_{k+2,l+2}} &
\CO_{\PP_S(E)/S}(k+2)\otimes \CB_{l+2} \ar[r]  &
\Coker \delta_{k+2,l+2} \to 0
}
$$
It follows from the diagram that $\sigma$ vanishes on $\Ker \delta_{k,l}$ and $\Coker \delta_{k,l}$,
hence both are supported on $\CX \subset \PP_S(E)$, the zero locus of $\sigma$.
On the other hand, the sheaf $\CO_{\PP_S(E)/S}(k-1)\otimes \CB_{l-1}$ is locally free,
hence has no torsion, thus $\Ker \delta_{k,l} = 0$.

Further, restricting this diagram to $\CX$ we deduce the desired exact sequences.
Finally, the periodicity isomorphisms follow immediately from the definition
of sheaves $\CE_{k,l}$ and $\CE'_{k,l}$.
\end{proof}

Gluing the short exact sequences of the lemma we obtain the following long exact sequences
of left and right $\CB_0$-modules on $\CX$:
\begin{eqnarray}
\dots \to
\CO_{\CX/S}(k-2)\otimes\CB_{l-2} \to
\CO_{\CX/S}(k-1)\otimes\CB_{l-1} \to
\CO_{\CX/S}(k)\otimes\CB_l \to
\CE_{k,l} \to 0,\label{e}\\
\dots \to
\CO_{\CX/S}(k-2)\otimes\CB_{l-2} \to
\CO_{\CX/S}(k-1)\otimes\CB_{l-1} \to
\CO_{\CX/S}(k)\otimes\CB_{l} \to
\CE'_{k,l} \to 0.\label{ep}
\end{eqnarray}

On the other hand, we have the following resolutions of $i_*\CE_{k,l}$ and $i_*\CE'_{k,l}$ on $\PP_S(E)$:
\begin{equation}\label{iseep}
\begin{array}{l}
0 \to \CO_{\PP_S(E)/S}(k-1)\otimes \CB_{l-1} \xymatrix@1{\ar[r]^{\delta_{k,l}} & }
\CO_{\PP_S(E)/S}(k)\otimes\CB_{l} \to
i_*\CE_{k,l} \to 0,\\
0 \to \CO_{\PP_S(E)/S}(k-1)\otimes \CB_{l-1} \xymatrix@1{\ar[r]^{\delta'_{k,l}} & }
\CO_{\PP_S(E)/S}(k)\otimes\CB_{l} \to
i_*\CE'_{k,l} \to 0.
\end{array}
\end{equation}

\begin{lemma}\label{eep}
The sheaves $\CE_{k,l}$ and $\CE'_{k,l}$ on $(\CX,\CB_0)$ and $(\CX,\CB_0^\opp)$ are
locally projective over $(S,\CB_0)$ and flat over $(S,\CB_0^\opp)$ respectively.
Moreover, there is an isomorphism of functors $\D^b(S,\CB_0) \to \D^b(\CX)$.
$$
\RCHom_{\CB_0}(\CE_{k,l},p^*F) \cong p^*F\otimes_{\CB_0}\CE'_{-k-1,-l-1}.
$$
\end{lemma}
\begin{proof}
Let $F$ be a sheaf of $\CB_0$-modules on $S$. Since $p:\CX\to S$ is flat,
$p^*F$ is a pure sheaf of $\CB_0$-modules on $\CX$. Further,
using periodicity isomorphism~(\ref{period}), the formula
for the twisted pullback functor of a regular closed embedding
$$
i^!(-) \cong
i^*(-)\otimes\CN_{\CX/\PP_S(E)}[-1] \cong
i^*(-)\otimes(\CO_{\CX/S}(2)\otimes \CL^*)[-1],
$$
and the duality theorem (see \cite{H}), we deduce
\begin{multline*}
i_*\RCHom_{\CB_0}(\CE_{k,l},p^*F) \cong
i_*\RCHom_{\CB_0}(\CE_{k+2,l+2},p^*F\otimes(\CO_{\CX/S}(2)\otimes \CL^*)) \cong
\\ \cong
i_*\RCHom_{\CB_0}(\CE_{k+2,l+2},i^!\pi^*F[1]) \cong
\RCHom_{\CB_0}(i_*\CE_{k+2,l+2},\pi^*F)[1].
\end{multline*}
Applying the functor $\RCHom_{\CB_0}(-,\pi^*F)$ to the resolution~(\ref{iseep}) of $i_*\CE_{k+2,l+2}$
and taking into account lemma~\ref{rhomb} we obtain the following distinguished triangle
$$
\pi^*F \otimes_{\CB_0} (\CO_{\PP_S(E)/S}(-k-2)\otimes\CB_{-l-2}) \stackrel{\delta'}\to
\pi^*F \otimes_{\CB_0} (\CO_{\PP_S(E)/S}(-k-1)\otimes\CB_{-l-1}) \to
\RCHom_{\CB_0}(i_*\CE_{k+2,l+2},\pi^*F)[1].
$$
On the other hand, the arguments of the lemma~\ref{kcd} show that the map $\delta'$ here
is an embedding, because the sheaf $\pi^*F$ has no torsion on $\CX$.
Since $i_*$ is exact and conservative, it follows that
$\RCHom_{\CB_0}(\CE_{k,l},p^*F)$ is a pure sheaf.
Since this is true for any $F \in \Coh(S,\CB_0)$, the sheaf $\CE_{k,l}$ is
locally projective over $(S,\CB_0)$.

Similarly, we have
$$
i_*(p^*F\otimes_{\CB_0}\CE'_{s,t}) \cong
i_*(i^*\pi^*F\otimes_{\CB_0}\CE'_{s,t}) \cong
\pi^*F \otimes_{\CB_0} i_*\CE'_{s,t}.
$$
Applying the functor $\pi^*F\otimes_{\CB_0}-$ to the resolution~(\ref{iseep}) of $i_*\CE'_{s,t}$ we deduce
the following distinguished triangle
$$
\pi^*F \otimes_{\CB_0} (\CO_{\PP_S(E)/S}(s-1)\otimes\CB_{t-1}) \stackrel{\delta'}\to
\pi^*F \otimes_{\CB_0} (\CO_{\PP_S(E)/S}(s)\otimes\CB_{t}) \to
\pi^*F \otimes_{\CB_0} i_*\CE'_{s,t}.
$$
By the same arguments as above we deduce that $p^*F\otimes_{\CB_0}\CE'_{s,t}$
is a pure sheaf, hence $\CE'_{s,t}$ is flat over $(S,\CB_0^\opp)$.

Further, comparing the above resolutions of $i_*\RCHom_{\CB_0}(\CE_{k,l},p^*F)$
and $i_*(p^*F\otimes_{\CB_0}\CE'_{s,t})$ it is easy to see that they coincide
for $(s,t) = (-k-1,-l-1)$. Therefore, we have a functorial isomorphism
of resolutions which gives a functorial isomorphism
$i_*\RCHom_{\CB_0}(\CE_{k,l},p^*F) \cong i_*(p^*F\otimes_{\CB_0}\CE'_{s,t})$,
(we use here the vanishing
$\Hom(\pi^*F \otimes_{\CB_0} (\CO_{\PP_S(E)/S}(-k-1)\otimes\CB_{-l-1}),
\pi^*F \otimes_{\CB_0} (\CO_{\PP_S(E)/S}(-k-2)\otimes\CB_{-l-2})) = 0$
which follows easily from lemma~\ref{easy}).
Since $i_*$ is conservative we obtain a functorial in $F$ isomorphism
$$
\RCHom_{\CB_0}(\CE_{k,l},p^*F) \cong p^*F\otimes_{\CB_0}\CE'_{s,t},
$$
for any $F \in \Coh(S,\CB_0)$. A standard argument then shows
that the above functorial isomorphism holds for any $F \in \D^b(S,\CB_0)$.
\end{proof}

\begin{lemma}\label{elf}
The sheaves $\CE_{k,l}$ and $\CE'_{k,l}$ are locally free over $\CO_\CX$ of rank $2^{n-2}$.
\end{lemma}
\begin{proof}
Let $x\in\CX$ be a point and $s = p(x)$. Let $V$ be the fiber of $E$ at $s$
and $v\in V$ a vector corresponding to the point $i(x) \in \PP(V)$.
Let $q = \sigma_s \in S^2V^*$ and consider the Clifford algebra $B_q$.
It follows from exact sequences~(\ref{e}) and (\ref{ep}) and from the definition
of $\delta$ that the fibers of $\CE_{k,l}$ and $\CE'_{k,l}$ at a point $x\in\CX$
coincide with the cokernel of the left and right multiplication maps
$B_{q0} \stackrel{v}\to B_{q1}$ (if $l$ is odd) or
$B_{q1} \stackrel{v}\to B_{q0}$ (if $l$ is even).
But the rank of these maps equals $2^{n-2}$ by lemma~\ref{multv} since $q(v) = 0$.
\end{proof}

Now we consider the following functors
$$
\begin{array}{ll}
\Phi_{k,l}:\D^b(S,\CB_0) \to \D^b(\CX),
\qquad&
F \mapsto p^*F \otimes_{\CB_0} \CE'_{k,l},\\
\Psi_{k,l}:\D^b(\CX) \to \D^b(S,\CB_0),
\qquad&
G \mapsto p_*(G \otimes_{\CO_\CX} \CE_{k,l}\otimes_{\CO_\CX}\det E^*[n-2]).
\end{array}
$$

Note that both functors preserve the bounded derived category by lemma~\ref{eep}
and lemma~\ref{elf}.

\begin{remark}\label{phiker}
Actually, the functors $\Phi_{k,l}$ and $\Psi_{k,l}$ are kernel functors
whose kernels are given by the sheaves $\CE'_{k,l}$ and $\CE_{k,l}\otimes_{\CO_\CX}\det E^*[n-2]$
on $(\CX,\CB_0) = (S,\CB_0)\times_S \CX$ respectively.
\end{remark}

\begin{proposition}\label{phieff}
For any $k,l$ the functor $\Phi_{k,l}:\D^b(S,\CB_0) \to \D^b(\CX)$ is fully faithful.
\end{proposition}

To prove the proposition we will check that the functor $\Psi_{1-k-n,1-l}$ is
left adjoint to $\Phi_{k,l}$ (lemma~\ref{ladj}), show that the composition
$\Psi_{1-k-n,1-l}\circ\Phi_{k,l}:\D^b(S,\CB_0) \to \D^b(S,\CB_0)$ is given by tensoring
with a $\CB_0$-bimodule $p_*(\CE'_{k,l}\otimes_{\CO_\CX}\CE_{1-k-n,1-l}\otimes_{\CO_\CX}\det E^*[n-2])$
(lemma~\ref{conv}) and check that this bimodule is isomorphic to $\CB_0$ (lemma~\ref{convepe}).


\begin{lemma}\label{ladj}
For any $k,l$ the functor $\Psi_{1-k-n,1-l} : \D^b(\CX) \to \D^b(S,\CB_0)$
is left adjoint to the functor $\Phi_{k,l}:\D^b(S,\CB_0) \to \D^b(\CX)$.
\end{lemma}
\begin{proof}
Take $F\in\D^b(S,\CB_0)$, $G\in\D^b(\CX)$.
Note that
$\Phi_{k,l}(F) =
p^*F \otimes_{\CB_0} \CE'_{k,l} \cong
\RCHom_{\CB_0}(\CE_{-k-1,-l-1},p^*F)$
by lemma~\ref{eep}.
Therefore
$\Hom(G,\Phi_{k,l}(F)) \cong \Hom(G,\RCHom_{\CB_0}(\CE_{-k-1,-l-1},p^*F))$.
Now we are going to use the duality theorem to rewrite this.
For this we replace the pullback functor $p^*$
with the twisted pullback functor $p^!$ using
$p^!F \cong p^*F \otimes \omega_{\CX/S}[n-2]$
and
$\omega_{\CX/S} \cong
\omega_{\PP_S(E)/S} \otimes \CO_{\CX/S}(2)\otimes \CL^* \cong
\CO_{\CX/S}(2-n)\otimes\det E^*\otimes \CL^*$.
Taking also into account the periodicity isomorphisms~(\ref{period}) we get
\begin{multline*}
\RCHom_{\CB_0}(\CE_{-k-1,-l-1},p^*F) \cong
\\ \cong
\RCHom_{\CB_0}(\CE_{-k-1,-l-1}\otimes \CO_{\CX/S}(2-n)\otimes\det E^*\otimes \CL^*[n-2],p^!F) \cong
\\ \cong
\RCHom_{\CB_0}(\CE_{1-k-n,1-l}\otimes\det E^*[n-2],p^!F).
\end{multline*}
Finally, applying the duality theorem we deduce
\begin{multline*}
\Hom(G,\RCHom_{\CB_0}(\CE_{1-k-n,1-l}\otimes\det E^*[n-2],p^!F)) \cong
\\ \cong
\Hom_{\CB_0}(G\otimes\CE_{1-k-n,1-l}\otimes\det E^*[n-2],p^!F) \cong
\\ \cong
\Hom_{\CB_0}(p_*(G\otimes\CE_{1-k-n,1-l}\otimes\det E^*[n-2]),F),
\end{multline*}
and an isomorphism
$\Hom(G,\Phi_{k,l}(F)) \cong \Hom(\Psi_{1-k-n,1-l}(G),F)$
follows.
\end{proof}

\begin{lemma}\label{conv}
The composition of functors
$\Psi_{1-k-n,1-l}\circ\Phi_{k,l}:\D^b(S,\CB_0) \to \D^b(S,\CB_0)$
is given by tensoring over $\CB_0$ by the complex of $\CB_0$-bimodules
$p_*(\CE'_{k,l}\otimes_{\CO_\CX}\CE_{1-k-n,1-l}\otimes_{\CO_\CX}\det E^*[n-2])$.
\end{lemma}
\begin{proof}
We have
\begin{multline*}
\Psi_{1-k-n,1-l}\circ\Phi_{k,l}(F) =
p_*((p^*F \otimes_{\CB_0} \CE'_{k,l}) \otimes_{\CO_X}\CE_{1-k-n,1-l}\otimes_{\CO_X}\det E^*[n-2]) \cong
\\ \cong
p_*(p^*F \otimes_{\CB_0} (\CE'_{k,l} \otimes_{\CO_X}\CE_{1-k-n,1-l}\otimes_{\CO_X}\det E^*[n-2])) \cong
\\ \cong
F \otimes_{\CB_0} p_*(\CE'_{k,l} \otimes_{\CO_X}\CE_{1-k-n,1-l}\otimes_{\CO_X}\det E^*[n-2]))
\end{multline*}
by associativity of the tensor product and the projection formula.
\end{proof}

\begin{lemma}\label{phiepok}
We have
$$
\Psi_{1-k-n,1-l}(\CO_{\CX/S}(k-t)) \cong
\CR^t_\FB \otimes_{\CB_0}\CB_{-l}
\quad\text{for all $t\ge 2-n$}.
$$
In particular, $\Psi_{1-k-n,1-l}(\CO_{\CX/S}(k-t)) = 0$ for $2-n\le t\le -1$.
%
%
%
\end{lemma}
\begin{proof}
Note that
$\Psi_{1-k-n,1-l}(\CO_{\CX/S}(k-t)) = p_*(\CO_{\CX/S}(k-t)\otimes_{\CO_\CX}\CE_{1-k-n,1-l}\otimes_{\CO_\CX}\det E^*[n-2])$
by definition of $\Psi_{1-k-n,1-l}$.
Tensoring resolution~(\ref{e}) by $\CO_{\CX/S}(k-t)\otimes_{\CO_\CX}\det E^*$
and taking into account lemma~\ref{poxs}
we obtain the following resolution
$$
\dots \to
\FA_{t+3}^*\otimes\CB_{-3-l} \to
\FA_{t+2}^*\otimes\CB_{-2-l} \to
\FA_{t+1}^*\otimes\CB_{-1-l} \to
p_*(\CO_{\CX/S}(k-t)\otimes_{\CO_\CX}\CE_{1-k-n,1-l}\otimes_{\CO_\CX}\det E^*[n-2])
$$
which coincides with the resolution~(\ref{crbl}) tensored over $\CB_0$ by $\CB_{-l}$.
\end{proof}

\begin{lemma}\label{convepe}
We have
$p_*(\CE'_{k,l}\otimes_{\CO_\CX}\CE_{1-k-n,1-l}\otimes_{\CO_\CX}\det E^*[n-2]) \cong \CB_0$.
\end{lemma}
\begin{proof}
Using the resolution~(\ref{ep}) and lemma~\ref{phiepok} above
we deduce that the LHS is quasiisomorphic to the complex
of sheaves of $\CB$-bimodules
$$
\dots \to
\CB_{l-2}\otimes_{\CO_\CX}(\CR^2_\FB\otimes_{\CB_0}\CB_{-l}) \to
\CB_{l-1}\otimes_{\CO_\CX}(\CR^1_\FB\otimes_{\CB_0}\CB_{-l}) \to
\CB_l\otimes_{\CO_\CX}(\CR^0_\FB\otimes_{\CB_0}\CB_{-l}) \to
0.
$$
which by~(\ref{resdb}) is  quasiisomorphic
to the $\CB_0$-bimodule $\CB_l\otimes_{\CB_0}\CB_0\otimes_{\CB_0}\CB_{-l}$,
which is isomorphic to $\CB_0$.
\end{proof}

Combining lemma~\ref{conv} and~\ref{convepe} we see that
$\Psi_{1-k-n,1-l}(\Phi_{k,l}(F)) \cong  F\otimes_{\CB_0}\CB_0 \cong F$,
hence the composition $\Psi_{1-k-n,1-l}\circ\Phi_{k,l}$ is isomorphic to the identity functor.
On the other hand, we have shown in lemma~\ref{ladj} that $\Psi_{1-k-n,1-l}$
is left adjoint to $\Phi_{k,l}$, hence $\Phi_{k,l}$ is fully faithful indeed.

To conclude the proof of theorem~\ref{dbcx} it remains to check that the category
$\Phi_{k,l}(\D^b(S,\CB_0))$ is right orthogonal to the collection~(\ref{rm2col})
and that together they generate $\D^b(\CX)$.

%

\begin{lemma}\label{phiro}
The subcategory $\Phi_{k,l}(D^b(S,\CB_0))$ is right orthogonal to the subcategory $(\ref{rm2col})$.
\end{lemma}
\begin{proof}
Recalling the definition of the functor $\Psi_{1-k-n,1-l}$ and using
the projection formula we deduce that
$$
\Psi_{1-k-n,1-l}(\CO_{\CX/S}(k-t)\otimes p^*F) \cong
\Psi_{1-k-n,1-l}(\CO_{\CX/S}(k-t))\otimes p^*F.
$$
But $\Psi_{1-k-n,1-l}(\CO_{\CX/S}(k-t)) = 0$ for $2-n\le t\le -1$
by lemma~\ref{phiepok}. Therefore
$$
\Psi_{1-k-n,1-l}(\langle \CO_{\CX/S}(k+1),\dots,\CO_{\CX/S}(k+n-2) \rangle\boxtimes\D^b(S)) = 0.
$$
Recalling lemma~\ref{ladj} we deduce by adjunction that
$\Phi_{k,l}(\D^b(S,\CB_0))$ is right orthogonal to $(\ref{rm2col})$.
\end{proof}


\begin{lemma}\label{phicpsi}
The composition of functors
$\Phi_{k,l}\circ\Psi_{1-k-n,1-l}:\D^b(\CX) \to \D^b(\CX)$
is isomorphic to the functor
$G\mapsto p_{2*}(p_1^*G\otimes_{\CO}(p_1^*(\CE_{1-k-n,1-l}\otimes_{\CO}\det E^*[n-2])\otimes_{\CB_0}p_2^*\CE'_{k,l}))$,
where $p_1,p_2:\CX\times_S\CX \to \CX$ are the projections to the factors.
$$
\xymatrix@!=10pt{
\CX\times_S\CX \ar[d]_{p_1} \ar[rr]^{p_2} && \CX \ar[d]_{p} \\
\CX \ar[rr]^{p} && S
}
$$
\end{lemma}
\begin{proof}
Indeed,
$\Phi_{k,l}\circ\Psi_{1-k-n,1-l}(G) =
p^*p_*(G \otimes_{\CO} \CE_{1-k-n,1-l}\otimes_{\CO}\det E^*[n-2]) \otimes_{\CB_0} \CE'_{k,l}$.
By the flat base change we have $p^*p_* \cong p_{2*}p_1^*$, and applying the projection formula
for $p_2$ we deduce the claim.
\end{proof}

Let us denote by $K$ the truncation of the resolution of the diagonal~(\ref{resda})
on $\CX\times_S\CX$ in degree $n-3$,
so that we get the following left and right resolutions for $K$:
\begin{eqnarray}
K &=& \left\{
0 \to
\CR^{n-3}_\FA\boxtimes\CO_{\CX/S}(3-n) \to
\dots \to
\CR^{1}_\FA\boxtimes\CO_{\CX/S}(-1) \to
\CR^{0}_\FA\boxtimes\CO_\CX \to
\Delta_*\CO_\CX \to 0
\right\},\label{rightk}\\
K &=& \left\{
\dots \to
\CR^{n}_\FA\boxtimes\CO_{\CX/S}(-n) \to
\CR^{n-1}_\FA\boxtimes\CO_{\CX/S}(1-n) \to
\CR^{n-2}_\FA\boxtimes\CO_{\CX/S}(2-n) \to
0\label{leftk}
\right\}
\end{eqnarray}

\begin{lemma}\label{convk}
We have
$p_1^*(\CE_{1-k-n,1-l}\otimes_{\CO}\det E^*[n-2])\otimes_{\CB_0}p_2^*\CE'_{k,l} \cong
K\otimes (\CO_\CX(2-k-n)\boxtimes\CO_\CX(k+n-2))$.
\end{lemma}
\begin{proof}
Tensoring over $\CB_0$
the pullback via $p_1$ of the resolution~(\ref{e}) of $\CE_{1-k-n,1-l}$ twisted by $\det E^*$
and the pullback via $p_2$ of the resolution~(\ref{ep}) of $\CE'_{k,l}$
we deduce that
$p_1^*(\CE_{1-k-n,1-l}\otimes_{\CO}\det E^*[n-2])\otimes_{\CB_0}p_2^*\CE'_{k,l}$
is quasiisomorphic to the total complex of the following bicomplex:
$$
\left\{\vcenter{\xymatrix{
\dots \ar[r] &
\CB_{-1}(-1-k-n,k) \ar[r] &
\CB_{0}(-k-n,k) \ar[r] &
\CB_{1}(1-k-n,k) \\
\dots \ar[r] &
\CB_{-2}(-1-k-n,k-1) \ar[r] \ar[u] &
\CB_{-1}(-k-n,k-1) \ar[r] \ar[u] &
\CB_{0}(1-k-n,k-1) \ar[u] \\
\dots \ar[r] &
\CB_{-3}(-1-k-n,k-2) \ar[r] \ar[u] &
\CB_{-2}(-k-n,k-2) \ar[r] \ar[u] &
\CB_{-1}(1-k-n,k-2) \ar[u] \\
& \vdots \ar[u] & \vdots \ar[u] & \vdots \ar[u]
}}\right\}\otimes\det E^*
$$
But comparing~(\ref{gsok}) with formulas
(\ref{cb0}), (\ref{cb1}) and (\ref{fb_k}) it is easy to see that
$$
\CB_k\otimes\det E^* \cong \FB_{n-k}^*
\qquad\text{for all $k\le 1$}.
$$
Therefore, this bicomplex can be rewritten as
$$
\left\{\vcenter{\xymatrix{
\dots \ar[r] &
\FB^*_{n+1}(-1-k-n,k) \ar[r] &
\FB^*_{n}(-k-n,k) \ar[r] &
\FB^*_{n-1}(1-k-n,k) \\
\dots \ar[r] &
\FB^*_{n+2}(-1-k-n,k-1) \ar[r] \ar[u] &
\FB^*_{n+1}(-k-n,k-1) \ar[r] \ar[u] &
\FB^*_{n}(1-k-n,k-1) \ar[u] \\
\dots \ar[r] &
\FB^*_{n+3}(-1-k-n,k-2) \ar[r] \ar[u] &
\FB^*_{n+2}(-k-n,k-2) \ar[r] \ar[u] &
\FB^*_{n+1}(1-k-n,k-2) \ar[u] \\
& \vdots \ar[u] & \vdots \ar[u] & \vdots \ar[u]
}}\right\}
$$
Comparing rows of this bicomplex with resolutions~(\ref{cral})
we see that its total complex is quasiisomorphic to
$$
\left\{
\dots \to
\CR_\FA^{n}(2-k-n)\boxtimes\CO_\CX(k-2) \to
\CR_\FA^{n-1}(2-k-n)\boxtimes\CO_\CX(k-1) \to
\CR_\FA^{n-2}(2-k-n)\boxtimes\CO_\CX(k)
\to 0
\right\}
$$
But looking at~(\ref{leftk}) we see that this complex coincides
with $K\otimes (\CO_\CX(2-k-n)\boxtimes\CO_\CX(k+n-2))$.
\end{proof}

Now we can give a proof of theorem.

\begin{proof}
By proposition~\ref{phieff} the functor $\Phi_{k,l}$ is fully faithful and
by lemma~\ref{easy} and lemma~\ref{phiro} the collection
$$
\lan
\D^b(S,\CB_0),
p^*(\D^b(S))\otimes\CO_{\CX/S}(1),
p^*(\D^b(S))\otimes\CO_{\CX/S}(2),
\dots,
p^*(\D^b(S))\otimes\CO_{\CX/S}(n-2)\ran
$$
is semiorthogonal in $\D^b(\CX)$.
It remains to check that $\D^b(\CX)$ is generated by this collection.
Tensoring the right resolution~(\ref{rightk}) of $K$ with
$(\CO_\CX(2-k-n)\otimes\CO_\CX(k+n-2))$ and taking into account lemma~\ref{convk}
we see that
\begin{multline*}
p_1^*(\CE_{1-k-n,1-l}\otimes_{\CO}\det E^*[n-2])\otimes_{\CB_0}p_2^*\CE'_{k,l} \cong
\\ \cong
\left\{
0 \to
\CR_\FA^{n-3}(2-k-n)\boxtimes\CO_\CX(k+1) \to
\dots \to
\CR_\FA^{0}(2-k-n)\boxtimes\CO_\CX(k+n-2) \to
\Delta_*\CO_\CX \to 0
\right\}.
\end{multline*}

Take any $G \in \D^b(\CX)$. Note that
$$
\begin{array}{lll}
p_{2*}(p_1^*G \otimes \Delta_*\CO_\CX) &\cong& G,\\
p_{2*}(p_1^*G \otimes \CR_\FA^{s}(2-k-n)\boxtimes\CO_\CX(k+n-2-s)) &\cong&
p^*p_*(G \otimes \CR_\FA^{s}(2-k-n))\otimes\CO_\CX(k+n-2-s),
\end{array}
$$
and on the other hand, by lemma~\ref{phicpsi} we have
$$
p_{2*}(p_1^*G \otimes p_1^*(\CE_{1-k-n,1-l}\otimes_{\CO}\det E^*[n-2])\otimes_{\CB_0}p_2^*\CE'_{k,l}) \cong
\Phi_{k,l}\circ\Psi_{1-k-n,1-l}(G).
$$
It follows that $G$ admits a filtration with quotients of the form
$$
p^*p_*(G \otimes \CR_\FA^{s}(2-k-n))\otimes\CO_\CX(k+n-2-s) \in p^*(\D^b(S))\otimes\CO_\CX(k+n-2-s),
\qquad
s=0,\dots,n-3,
$$
and $\Phi_{k,l}(\Psi_{1-k-n,1-l}(G)) \in \Phi_{k,l}(\D^b(S,\CB_0))$.
Therefore $\D^b(\CX)$ is indeed generated by the desired collection.
\end{proof}

\section{Homological Projective Duality for the double Veronese variety}

We start this section with a brief reminder of the notion of the Homological Projective Duality.
The general reference for this is~\cite{K2}.
Assume that $X$ is an algebraic variety with a line bundle $\CO_X(1)$ on $X$.

\begin{definition}[\cite{K2}]
A {\sf Lefschetz decomposition} of the derived category $\D^b(X)$ is
a semiorthogonal decomposition of $\D^b(X)$ of the form
\begin{equation}\label{dbx}
\D^b(X) = \langle \CA_0,\CA_1(1),\dots,\CA_{\ix-1}(\ix-1) \rangle,
\qquad
0 \subset \CA_{\ix-1} \subset \CA_{\ix-2} \subset \dots \subset \CA_1 \subset \CA_0 \subset \D^b(X),
\end{equation}
where $0 \subset \CA_{\ix-1} \subset \CA_{\ix-2} \subset \dots \subset \CA_1 \subset \CA_0 \subset \D^b(X)$
is a chain of admissible subcategories of $\D^b(X)$.
\end{definition}

Let $\fa_k$ denote the right orthogonal to $\CA_{k+1}$ in $\CA_k$, so that
we have a semiorthogonal decomposition
$$
\CA_k = \langle\fa_k,\fa_{k+1},\dots,\fa_{\ix-1}\rangle.
$$
The categories $\fa_0,\fa_1,\dots,\fa_{\ix-1}$ are called {\sf primitive}
categories of the Lefschetz decomposition~(\ref{dbx}).

Let $f:X \to \PP(V)$ be a morphism into a projective space such that
$f^*(\CO_{\PP(V)}(1)) \cong \CO_X(1)$ and let $\CX \subset X\times\PP(V^*)$
be the universal hyperplane section of $X$
(i.e. the canonical divisor of bidegree $(1,1)$ in $X\times\PP(V^*)$).

\begin{definition}[\cite{K2}]
An algebraic variety $Y$ with a projective morphism $g:Y\to\PP(V^*)$
is called {\sf Homologically Projectively Dual} to $f:X\to \PP(V)$
with respect to a Lefschetz decomposition~$(\ref{dbx})$, if
there exists an object $\CE\in\D^b(\CX\times_{\PP(V^*)} Y)$ such that
the kernel functor
$\Phi = \Phi_\CE:\D^b(Y) \to \D^b(\CX)$ is fully faithful
and gives the following semiorthogonal decomposition
\begin{equation}\label{dbx1}
\D^b(\CX) = \langle \Phi(\D^b(Y)),\CA_1(1)\boxtimes\D^b(\PP(V^*)),\dots,
\CA_{\ix-1}(\ix-1)\boxtimes\D^b(\PP(V^*))\rangle.
\end{equation}
\end{definition}

For every linear subspace $L \subset V^*$ we consider
the corresponding linear sections of $X$ and $Y$:
$$
X_L = X\times_{\PP(V)}\PP(L^\perp),
\qquad
Y_L = Y\times_{\PP(V^*)}\PP(L),
$$
where $L^\perp \subset V$ is the orthogonal subspace to $L\subset V^*$.
Let $N = \dim V$.

The main property of Homologically Projectively Dual varieties is the following

\begin{theorem}[\cite{K2}]\label{hp}
If\/ $Y$ is Homologically Projectively Dual to $X$ then

\noindent$(i)$
$Y$ is smooth and $\D^b(Y)$ admits a dual Lefschetz decomposition
$$
\D^b(Y) = \lan \CB_{\jx-1}(1-\jx),\dots,\CB_{1}(-1),\CB_0\ran,\qquad
0 \subset \CB_{\jx-1} \subset \dots \subset \CB_1 \subset \CB_0 \subset \D^b(Y)
$$
with the same set of primitive subcategories:
$\CB_k = \langle\fa_0,\dots,\fa_{N-k-2}\rangle$;

\noindent$(ii)$
for any linear subspace $L\subset V^*$, $\dim L = r$, such that we have
$\dim X_L = \dim X - \dim L$, and $\dim Y_L = \dim Y + \dim L - N$\/
there exist a triangulated category $\CC_L$ and semiorthogonal decompositions
$$
\begin{array}{lll}
\D^b(X_L) &=& \langle \CC_L,\CA_{r}(1),\dots,\CA_{\ix-1}(\ix-r)\rangle\\
\D^b(Y_L) &=& \langle \CB_{\jx-1}(N-r-\jx),\dots,\CB_{N-r}(-1),\CC_L\rangle.
\end{array}
$$
\end{theorem}

Now we are going to apply theorem~\ref{hp} in the following situation.
Let $W$ be a vector space, $\dim W = n$, and $X = \PP(W)$.
Consider the standard exceptional collection
$$
\D^b(X) = \lan \CO_X(-1),\CO_X,\CO_X(1),\dots,\CO_X(n-2)\ran
$$
on the projective space $X = \PP(W) \cong \PP^{n-1}$ and put
$$
\CA_0 = \CA_1 = \dots = \CA_{\ix-2} := \lan \CO_X(-1),\CO_X \ran,
\qquad
\CA_{\ix-1} := \begin{cases}
\lan \CO_X(-1),\CO_X \ran, & \text{if $n = 2\ix$   is even}\\
\lan \CO_X(-1)\ran,        & \text{if $n = 2\ix-1$ is odd}
\end{cases}
$$
where $\ix = \lceil n/2\rceil$.
Then
\begin{equation}\label{ldx}
\D^b(X) = \lan\CA_0,\CA_1\otimes\CO_{\PP(S^2W)}(1),\dots,\CA_{\ix-1}\otimes\CO_{\PP(S^2W)}(\ix-1)\ran
\end{equation}
is a Lefschetz decomposition of $\PP(W)$ with respect to the double Veronese embedding $f:\PP(W) \to \PP(S^2W)$
(note that $f^*(\CO_{\PP(S^2W)}(1)) \cong \CO_{\PP(W)}(2)$).

Let $\CX \subset \PP(W)\times\PP(S^2W^*)$ be the universal hyperplane
section of $\PP(W)$ with respect to the double Veronese embedding $f:\PP(W) \to \PP(S^2W)$.
In other words, $\CX$ is the {\sf universal quadric} (the fiber of $\CX$ over a point
of the space $\PP(S^2W^*)$ is the corresponding quadric in $\PP(W)$).
Thus $\CX \to \PP(S^2W^*)$ is a quadric fibration.
Consider the corresponding sheaf of even parts of Clifford algebras on $\PP(S^2W^*)$
$$
\CB_0 =
\CO_{\PP(S^2W^*)} \oplus
\Lambda^2W \otimes \CO_{\PP(S^2W^*)}(-1) \oplus
\Lambda^4W \otimes \CO_{\PP(S^2W^*)}(-2) \oplus
\dots,
$$
constructed in section~\ref{gfb}

\begin{theorem}\label{xyhpd}
The noncommutative algebraic variety $Y = (\PP(S^2W^*),\CB_0)$
is Homologically Projectively Dual to $\PP(W)$ with respect
to the Lefschetz decomposition~$(\ref{ldx})$.
\end{theorem}
\begin{proof}
Rewriting the semiorthogonal decomposition of $\D^b(\CX)$
given by theorem~\ref{dbcx} in our case we get
$$
\D^b(\CX) = \langle \Phi_{0,0}(\D^b(Y)),\CA_1(1)\boxtimes\D^b(\PP(S^2W^*)),\dots,
\CA_{\ix-1}(\ix-1)\boxtimes\D^b(\PP(S^2W^*))\rangle.
$$
It remains to note that the functor $\Phi_{0,0}$ is a kernel functor
with kernel supported on the fiber product $\CX\times_{\PP(S^2W^*)} Y = (\CX,\CB_0)$
by remark~\ref{phiker}.
\end{proof}

A direct computation shows that the dual Lefschetz decomposition
of $\D^b(\PP(S^2W^*),\CB_0)$ in this case takes form
$$
\D^b(\PP(S^2W^*),\CB_0) =
\lan\CB_{-n^2},\CB_{1-n^2},\dots,\CB_{-1},\CB_0,\CB_1\ran,
$$
where
$$
\begin{array}{lcl}
\CB_1 &=&
W \otimes \CO_{\PP(S^2W^*)} \oplus
\Lambda^3W \otimes \CO_{\PP(S^2W^*)}(-1) \oplus
\Lambda^5W \otimes \CO_{\PP(S^2W^*)}(-2) \oplus
\dots,\\
\CB_{k+2} &=& \CB_k\otimes\CO_{\PP(S^2W^*)}(1).
\end{array}
$$
is the sequence of $\CB_0$-modules constructed in section~\ref{gfb}.

Now we are going to apply theorem~\ref{hp} to the Homologically Projectively Dual
varieties $X = \PP(W)$ and $Y = (\PP(S^2W^*),\CB_0)$.
Then for each vector subspace $L \subset S^2W^*$ the corresponding linear section
$X_L := \PP(W) \times_{\PP(S^2W^*)} \PP(L^\perp)$ is the intersection of quadrics
in $\PP(W)$ parameterized by the space $L$, and
$Y_L := (\PP(S^2W^*),\CB_0) \times_{\PP(S^2W^*)} \PP(L) = (\PP(L),\CB_0)$.
The subspace $L$ is admissible if and only if the intersection of quadrics $X_L$
is complete. Thus we deduce the following description of the derived category
of a complete intersection of quadrics:

\begin{theorem}\label{dbxl}
For any vector subspace $L\subset S^2W^*$ such that the corresponding
intersection of quadrics $X_L$ is complete there exists
a semiorthogonal decomposition
$$
\begin{array}{llll}
\D^b(X_L) &=&  \langle \D^b(\PP(L),\CB_0),\CO_{X_L}(1),\dots,\CO_{X_L}(n-2r)\rangle, &
\text{if $r = \dim L \le n/2$}\\
\D^b(\PP(L),\CB_0) &=&  \langle \CB_{n-2r},\dots,\CB_{-2},\CB_{-1},\D^b(X_L)\rangle, &
\text{if $r = \dim L \ge n/2$}
\end{array}
$$
and an equivalence $\D^b(X_L) \cong \D^b(\PP(L),\CB_0)$ if $\dim L = n/2$,
where $\D^b(\PP(L),\CB_0)$ is the derived category of sheaves of $\CB_0$-modules on $\PP(L)$.
\end{theorem}

In the case when $r = \dim L = 1$ the space $X_1 = X_L$ is a single quadric and $\D^b(\PP(L),\CB_0)$
is the derived category of $B_0$-modules where $B_0$ is the even part of the Clifford algebra
of the quadric. If the quadric is smooth then the algebra $B_0$ is semisimple and we recover
the Kapranov's description of the derived category of a smooth quadric.

\begin{corollary}
If $X_1$ is a quadric in $\PP^{n-1}$ then we have a semiorthogonal decomposition
$$
\D^b(X_1) = \lan \D^b(B_0),\CO_{X_1}(1),\dots,\CO_{X_1}(n-2)\ran,
$$
where $B_0$ is the even part of the Clifford algebra of the quadric,
and $\D^b(B_0)$ is the derived category of $B_0$-modules.
If $X_1$ is smooth then $\D^b(B_0)$ is generated by one {\rm(}for odd $n${\rm)}
or two {\rm(}for even $n${\rm)} exceptional objects.
\end{corollary}

In the case when $r = \dim L = 2$ the space $X_2 = X_L$ is an intersection of two quadrics.
If the intersection is complete and smooth we apply corollaries~\ref{dc} and \ref{sc} and obtain
the following description of the derived category $\D^b(X_L)$:

\begin{corollary}
If $X_2$ is a smooth complete intersection of two quadrics in $\PP^{n-1}$
then we have a semior\-thogonal decomposition
$$
\D^b(X_1) = \lan \D^b(C),\CO_{X_1}(1),\dots,\CO_{X_1}(n-4)\ran,
$$
where $C$ is either a twofold covering of $\PP^1$ ramified at the critical values of $p$ {\rm(}for even $n${\rm)},
or $\PP^1$ with $\ZZ/2\ZZ$-stack structure at the the critical values of $p$ {\rm(}for odd $n${\rm)}.
\end{corollary}

\end{document}